\documentclass[review,12pt]{elsarticle}
\usepackage[margin=1in]{geometry}

\usepackage{mathptmx}
\usepackage{helvet}
\usepackage{courier}
\usepackage{type1cm}
\usepackage{makeidx}
\usepackage{graphicx}
\usepackage{multicol}
\usepackage[bottom]{footmisc}
\usepackage{amsmath}
\usepackage{amssymb,amscd}
\usepackage{amsfonts}
\usepackage{enumitem}
\usepackage[
colorlinks=false,
citecolor=black,
urlcolor=black,
linktocpage,
pagebackref
]{hyperref}
\usepackage{comment}
\usepackage{algorithmic}
\usepackage[ruled,vlined,linesnumbered]{algorithm2e}

\SetKwRepeat{Do}{do}{while}
\def\bx{{\bf x}}   %bold x
\def\cmb{{\omega}}   %weights for combination of metric
\def\lim{{\xi}}    %limiting function

\def\hr{$hr-$}
\def\hra{$hr-$adaptivity}
\def\ra{$r-$adaptivity}
\def\ha{$h-$adaptivity}
\def\pa{$p-$adaptivity}
\def\sh{\emph{shape}}
\def\sz{\emph{size}}
\def\shsz{\emph{shape}+\emph{size}}

\DeclareMathOperator*{\argmin}{arg\,min}
\journal{Journal}

\begin{document}

\begin{frontmatter}

\title{$hr-$adaptivity for nonconforming high-order meshes with the target
matrix optimization paradigm}

\author{Veselin Dobrev \fnref{llnl}}
\author{Patrick Knupp \fnref{dihedral}}
\author{Tzanio Kolev \fnref{llnl}}
\author{Ketan Mittal \corref{cor1}  \fnref{llnl}}
\author{Vladimir Tomov \fnref{llnl}}
\fntext[llnl]
{Lawrence Livermore National Laboratory, 7000 East Avenue, Livermore, CA 94550}
\fntext[dihedral]
{Dihedral LLC, Bozeman, MT 59715}
\cortext[cor1]
{Corresponding author, mittal3@llnl.gov}
\tnotetext[l_title]
{Performed under the auspices of the U.S. Department of Energy under
Contract DE-AC52-07NA27344 (LLNL-JRNL-814656)}

\address{}
\begin{abstract}
We present an \hra\ framework for optimization of high-order meshes.  This work
extends the $r$-adaptivity method for mesh optimization by Dobrev et al.
\cite{TMOP2020}, where we utilized the Target-Matrix Optimization Paradigm
(TMOP) to minimize a functional that depends on each element's current and
target geometric parameters: element {\em aspect-ratio}, {\em size}, {\em
skew}, and {\em orientation}.  Since fixed mesh topology limits the ability to
achieve the target size and aspect-ratio at each position, in this paper we
augment the $r$-adaptivity framework with nonconforming adaptive mesh
refinement to further reduce the error with respect to the target geometric
parameters.  The proposed formulation, referred to as \hra, introduces
TMOP-based quality estimators to satisfy the aspect-ratio-target via {\em
anisotropic} refinements and size-target via {\em isotropic} refinements in
each element of the mesh.  The methodology presented is purely algebraic,
extends to both simplices and hexahedra/quadrilaterals of any order, and
supports nonconforming isotropic and anisotropic refinements in 2D and 3D.
Using a problem with a known exact solution, we demonstrate the effectiveness
of \hra\ over both $r-$ and \ha\ in obtaining similar accuracy in the solution
with significantly fewer degrees of freedom.  We also present several examples
that show that \hra\ can help satisfy geometric targets even when \ra\ fails to
do so, due to the topology of the initial mesh.
\end{abstract}

\end{frontmatter}

%Three article highlights as bullets with each less than 120 characters.
\section*{Article Highlights}
\begin{itemize}
\item A novel \hra\ method is developed for optimizing nonconforming high-order
meshes.
\item \hra\ is more effective in satisfying target geometric parameters of the
mesh in comparison to $h-$ and \ra\ techniques.
\item \hra\ meshes require fewer degrees of freedom for a given solution
accuracy in comparison to $h$- and \ra\ meshes.
\end{itemize}
%=====================================================================
%=====================================================================
\section{Introduction}
Mesh quality impacts the fidelity and robustness of the solution of a partial
differential equation (PDE) obtained by methods such as the finite element
method (FEM), finite volume method (FVM), and the spectral element method
(SEM).  The quality of a mesh depends on geometrical parameters such as shape
and size of each element in the mesh. The ideal values of these parameters
depend on factors such as the PDE type, the domain on which the problem is
being solved, element types in the mesh (simplices or
quadrilaterals/hexahedra), etc.  There are three popular approaches for
improving the quality of a given mesh without fully reconstructing it, namely,
\ra, \ha, and \pa.

The \ra\ methods move the nodes of the mesh elements without changing the mesh
topology.  Laplacian smoothing, where each node is typically moved as a linear
function of the positions of its neighbors \cite{vollmer1999improved,
field1988laplacian, taubin2001linear}, optimization-based smoothing, where a
functional based on elements' geometrical parameters is minimized
\cite{Knupp2012, gargallo2015optimization, mittal2019mesh, dobrev2019target,
Greene2017, Peiro2018}, and equidistribution with respect to an appropriate
metric tensor \cite{Huang1994, Huang2010} are three popular approaches for \ra.
Methods based on \ra\ enable us to improve the mesh quality to better capture
the solution of the PDE of interest. In certain cases, they can also improve
the conditioning of the linear system resulting from discretizing the PDE of
interest on the given mesh (e.g., \cite{mittal2019mesh}).  However, there are
two key limitations of \ra. First, the effectiveness of \ra\ methods is
restricted by the topology of the given mesh and can fail to achieve the target
parameters in the region of interest. Second, once a mesh has been modified by
moving its nodal positions, the PDE solution has to be transferred from the
original mesh to the improved mesh to continue the PDE solve process. This
grid-to-grid transfer can be computationally expensive, especially for
high-order discretizations.

In contrast to \ra, $h-$ and $p-$adaptivity based methods do not move the nodal
position of the elements, but instead add resolution to the domain by locally
refining individual elements of the mesh. In \ha, an element in the mesh can be
split into more elements \cite{cerveny2019nonconforming}, while in \pa, the
polynomial order used to represent the solution
%(for Galerkin-based methods)
can be changed in each element to provide more resolution where needed
\cite{barros2004error}.  In each of these approaches, typically an error
estimator is employed to determine the elements that should be refined to
improve the accuracy of the solution.  The advantage of $h-$ and $p-$adaptivity
is that each element can be refined to obtain arbitrary resolution in any
region of the mesh, irrespective of the shape and size of the original mesh.
Additionally, once an element is refined, it is trivial to interpolate the
solution from the original element (\emph{master}) to its refined counterparts
(\emph{slaves}) using interpolation matrices. However, there are two key
limitations of $h-$ and $p-$adaptivity. First, these methods cannot control the
geometric shape of the mesh elements. Thus, if the error in the solution is due
to the skewness or orientation of the element (e.g., magnetohydrodynamic
applications), $h$- or $p$-adaptivity methods cannot be used to produce
arbitrary accuracy in the solution.  Second, as each element is refined to
increase the resolution in the domain, the computational cost of the
calculation increases.

Based on the advantages and disadvantages of each method, a combination of
$r-$adaptivity with $h-$ and/or $p-$adaptivity can help circumvent the
shortcomings of each of these individual methods.  In this work, we focus on
$hr-$adaptivity; $r\!p-$ and $hr\!p-$adaptivity will be explored in future
papers.  A survey of the literature shows that \hra\ methods have been
effective for mesh optimization with applications ranging from the
Schr{\"o}dinger equations \cite{mackenzie2020hr} to ocean modeling
\cite{piggott2005h}. These methods typically employ an error-based function
that is minimized by moving the nodal positions (\ra) followed by
refining/derefining each element (\ha) by elemental operations such as node
insertion, node removal, or edge swapping.  A non-exhaustive list of
publications and recent advances on the subject of \hra\ is given by
\cite{mackenzie2020hr, piggott2005h, jahandari2020finite,
piggott2009anisotropic, ong2013hr, mostaghimi2015anisotropic, piggott2008new,
walkley2002anisotropic, antonietti2020hr, edwards1993h}.  All of the existing
methods \hr\ methods are either restricted to low-order meshes (first- or
second-order), specific element types (simplices or quadrilaterals/hexahedra),
or specific $h-$refinement types (typically isotropic).

In this work, we present a high-order \hra\ method based on the Target-Matrix
Optimization Paradigm (TMOP).  The algorithm is applicable to curved meshes of
any order.  It supports isotropic refinements for simplices, and both isotropic
and anisotropic anisotropic refinements for quadrilaterals/hexahedra, leading
to nonconforming 2D and 3D meshes with hanging nodes.  This paper demonstrates
that TMOP can be used in \hra\ settings to adapt the shape and size of the mesh
elements, and ultimately improve the accuracy, compared to standalone \ha\ and
\ra, with which the solution of a PDE can be represented on the given mesh.

The rest of the paper is organized as follows. In Section \ref{sec_prelim} we
review TMOP for \ra\ \cite{TMOP2020} and the adaptive mesh refinement framework
(\ha) \cite{cerveny2019nonconforming} that are the starting points of this
work. In Section \ref{sec_method} we describe how we integrate the \ha\
framework with TMOP to construct a novel \hra\ technique. Finally, we present
several numerical experiments in Section \ref{sec_apps} that show the
effectiveness of our framework in improving the computational efficiency and
accuracy of the numerical solution in comparison to \ha\ and \ra.  Conclusions
are presented in Section \ref{sec_conc}.

%=====================================================================
%=====================================================================
\section{Overview of TMOP, \ra\ and \ha} \label{sec_prelim}

This section provides a basic description of the \ra\ and \ha\ algorithms
that are the starting point of this work. We only focus on the aspects
that are relevant to the final \hra\ method.

%==========================

\subsection{Discrete mesh representation} \label{subsec_mesh}

In our finite element based framework, the domain $\Omega \in \mathbb{R}^d$ is
discretized as a union of curved mesh elements of order $k$.  To obtain a
discrete representation of these elements, we select a set of scalar basis
functions $\{ \bar{w}_i \}_{i=1}^{N_w}$ on the reference element $\bar{E}$.
This basis spans the space of all polynomials of degree at most $k$ on the
given element type (quadrilateral, tetrahedron, etc.).  The position of an
element $E$ in the mesh $\mathcal{M}$ is fully described by a matrix
$\mathbf{x}_E$ of size $d \times N_w$ whose columns represent the coordinates
of the element control points (a.k.a. nodes or element degrees of freedom).
Given $\mathbf{x}_E$, we introduce the map $\Phi_E:\bar{E} \to \mathbb{R}^d$
whose image is the element $E$:
\begin{equation}
\label{eq_x}
x(\bar{x}) =
   \Phi_E(\bar{x}) \equiv
   \sum_{i=1}^{N_w} \mathbf{x}_{E,i} \bar{w}_i(\bar{x})\,,
   \qquad \bar{x} \in \bar{E},
\end{equation}
where we used $\mathbf{x}_{E,i}$ to denote the $i$-th column of $\mathbf{x}_E$,
i.e., the $i-$th degree of freedom of element $E$.  To ensure continuity
between mesh elements, we define a global vector $\mathbf{x}$ of mesh positions
that contains the $\mathbf{x}_E$ control points for every element.

%Since the mesh positions $\bx \in \mathcal{V}$, mesh continuity is ensured
%using a gather-scatter operator (e.g., $P$ in Section 5.3 of
%\cite{MFEMPaper2019}) that maps the global degrees of freedom ($\bx$) to local
%element-by-element degrees of freedom ($\mathbf{x}_{E}$).

%==========================

\subsection{TMOP for \ra} \label{subsec_ra}
The objective of the \ra\ process is to optimize the mesh using information
from a discrete function, e.g., a finite element solution function that is
defined with respect to the initial mesh.  In TMOP, \ra\ is achieved by
incorporating the discrete data into the target geometrical configuration.  In
this subsection we summarize the main components of the TMOP approach; all
details of the specific method we use are provided in \cite{TMOP2020}.

The major concept of TMOP is the user-specified transformation ($W$), from
reference-space to target coordinates, which represents the ideal geometric
properties of every mesh point.  The construction of this transformation is
guided by the fact that any Jacobian matrix can be written as a combination of
four components:
\begin{equation}
\label{eq_W}
W = \underbrace{\zeta}_{\text{[volume]}} \underbrace{R}_{\text{[rotation]}}
\underbrace{Q}_{\text{[skewness]}} \underbrace{D}_{\text{[aspect-ratio]}}.
\end{equation}
Further details about this decomposition and a thorough discussion on how
TMOP's target construction methods encode geometric information into the target
matrix $W$ is given by Knupp in \cite{knupp2019target}.  For \ra\ in PDE-based
applications, the geometric parameters of \eqref{eq_W} are typically
constructed as discrete functions using the discrete solution available on the
initial mesh.  As the nodal coordinates change during \ra, the discrete
functions have to be mapped from the original mesh to the updated mesh to
ensure that $W$ can be constructed at each reference point, see Section 4 in
\cite{IMR2018}.

Using \eqref{eq_x}, the Jacobian of the mapping $\Phi_E$ at any reference point
$\bar{\bx} \in \bar{E}$ from the reference-space coordinates to the current
physical-space coordinates is defined as
\begin{equation}
\label{eq_A}
  A(\bar{\bx}) = \frac{\partial \Phi_E}{\partial \bar{\bx}} =
    \sum_{i=1}^{N_w} \mathbf{\bx}_{E,i} [ \nabla \bar{w}_i(\bar{\bx}) ]^T \,.
\end{equation}
\begin{figure}[tb!]
\centerline{
  \includegraphics[width=0.5\textwidth]{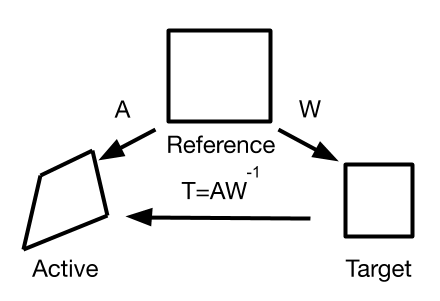}}
\caption{Schematic representation of the major TMOP matrices.}
\label{fig_tmop}
\end{figure}
Combining \eqref{eq_A} and \eqref{eq_W}, the transformation from the
target coordinates to the current physical coordinates (see Fig. \ref{fig_tmop})
is defined as
\begin{equation}
\label{eq_T}
T = AW^{-1}.
\end{equation}
With the transformation $T$, a quality metric $\mu(T)$ is used to measure the
difference between $A$ and $W$ in terms of the geometric parameters of interest
specified in \eqref{eq_W}. For example, $\mu_2=\mid T \mid^2/2\tau-1$ is a
\emph{shape} metric that depends on the skew and aspect ratio components, but
is invariant to orientation and volume.  Here, $\mid T \mid$ is the Frobenius
norm of $T$ and $\tau=\text{det}(T)$.  Similarly, $\mu_{55}=(\tau-1)^2$ is a
\emph{size} metric that depends only on the volume of the element. We also have
\emph{shape}$+$\emph{size} metrics such as $\mu_7 = \mid T-T^{-t}\mid^2$ and
$\mu_9=\tau\mid T-T^{-t}\mid^2$ that depend on volume, skew and aspect ratio,
but are invariant to orientation.

The quality metric $\mu(T)$ is used for \ra\ by minimizing the TMOP objective
function: \begin{equation}
\label{eq_F_full}
  F(\bx) = \sum_{E \in \mathcal{M}}F_E(\bx_E) = \sum_{E(\bx_E)} \int_{E_t}
  \cmb(\bx) \mu_{i}(T(\bx)) d\bx_t,
\end{equation}
where $F$ is a sum of the TMOP objective function for each element in the mesh ($F_E$),
$\cmb$ is a user-prescribed spatial weight
and $E_t$ is the target element corresponding to the physical element $E$.
In \eqref{eq_F_full}, the integral is computed as
\begin{equation}
\label{eq_vm}
  \sum_{E \in \mathcal{M}} \int_{E_t} \cmb(\bx_t) \mu(T(\bx_t)) d\bx_t = \frac{1}{N_E}
  \sum_{E \in \mathcal{M}} \sum_{\bx_q \in E_t}
                           w_q\,\det(W(\bar{\bx}_q))\, \cmb(\bx_q) \mu(T(\bx_q)),
\end{equation}
where $\mathcal{M}$ is the current mesh with $N_E$ elements, $w_q$ are the
quadrature weights, and the point $x_q$ is the image of the reference
quadrature point $\bar{x}_q$ in the target element $E_t$.

\begin{figure}[t!]
\begin{center}
$\begin{array}{ccc}
\includegraphics[height=0.3\textwidth]{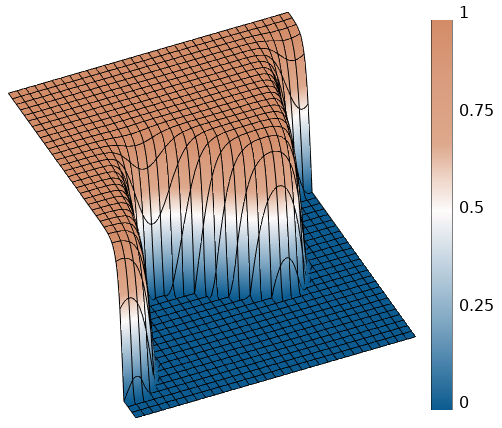} &
\includegraphics[height=0.3\textwidth]{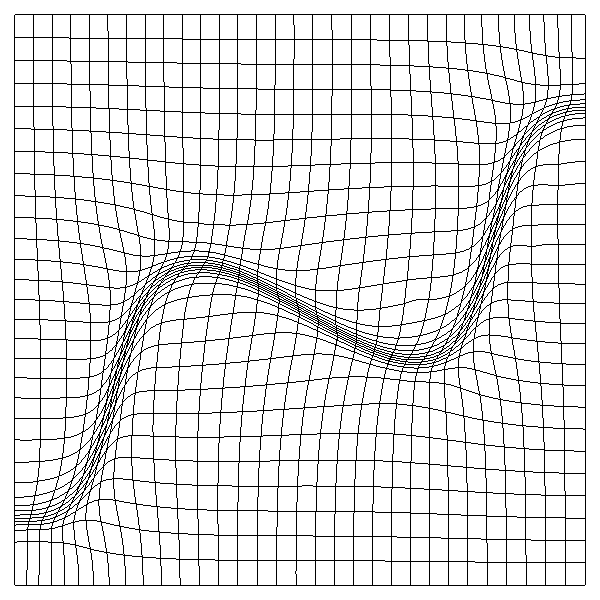} \\
\end{array}$
\end{center}
\vspace{-7mm}
\caption{(left) Material indicator function on the initial uniform mesh
         and (right) the optimized mesh.}
\label{fig_2matind}
\end{figure}

In practice, one may use a combination of multiple metrics with different
spatial weights to control different geometric parameters in different regions
of the mesh. In that case, the TMOP objective function is defined as a
combination of integrals associated with each metric \cite{TMOP2020}.  Using an
optimization problem solver (e.g., Newton's method and L-BFGS), optimum nodal
locations can be determined by minimizing the TMOP objective function
\eqref{eq_F_full}.  Figure \ref{fig_2matind} shows an example of \ra\ to a
discrete material indicator using TMOP to control the aspect-ratio and the size
of the elements in a mesh. In this example, the material indicator function is
used to define discrete functions for aspect-ratio and size targets in
\eqref{eq_W}, and the mesh is optimized using a \shsz\ metric.

%==========================

\subsection{Basics of \ha} \label{subsec_ha}
The \ha\ component of this work is based on the adaptive mesh refinement
framework by Cerveny et al. \cite{cerveny2019nonconforming}. Therein, the
authors present a highly scalable approach for unstructured nonconforming \ha\
that can be used for high-order curved meshes consisting of triangles,
quadrilaterals, tetrahedra, and hexahedra. The methods support the entire de
Rham sequence of finite element spaces, at arbitrarily high-order.

\begin{figure}[t!]
\centerline{
\includegraphics[width=0.80\textwidth]{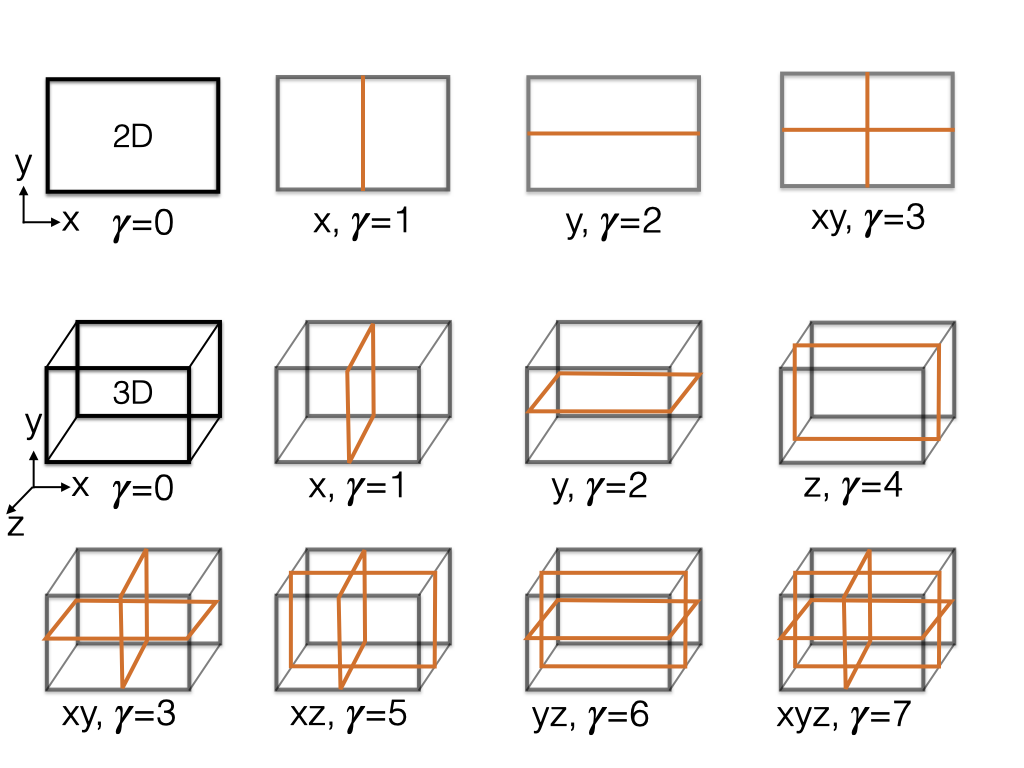} \hfil}
\caption{Different refinements options in 2D/3D.}
\label{fig_amr_refs}
\end{figure}

For the purposes of this work, we show all possible refinements for a
quadrilateral and a hexahedron element in Figure \ref{fig_amr_refs}, where we
have used $\gamma$ to represent the refinement type.  For a hexahedron, $\gamma
= 1\dots 6$ represent different types of anisotropic refinements and $\gamma =
7$ represents isotropic refinement.  The \ha\ framework by Cerveny et al. also
supports derefinement where elements introduced by refinement can be removed to
restore the parent element that they had originally emanated from. Note that
for derefinement, only the elements that share the same parent can be combined
together to coarsen the mesh, and a parent can be restored only by coarsening
all its children and not a subset thereof.  Depending on the specific target,
our \hra\ method can utilize all different $\gamma$ options for refinement and
derefinement.

\begin{figure}[t!]
\centerline{
  \includegraphics[width=0.80\textwidth]{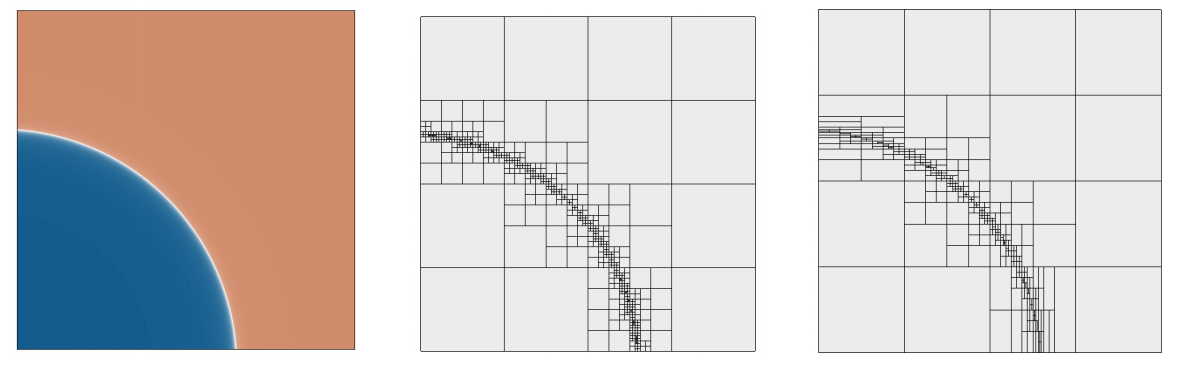} \hfil}
\caption{Discrete function in a domain (left) used for \ha\ with
isotropic refinements (center) and anisotropic refinements (right).}
\label{fig_amrh}
\end{figure}

In practice, an error estimator \cite{babuska1979reliable,
ainsworth2011posteriori} is used to determine the elements that should be
refined or derefined to improve the accuracy of the solution on a given mesh.
Figure \ref{fig_amrh} shows an example of \ha\ with isotropic and anisotropic
refinements based on the spatial gradients of the function.  Our \hra\ method
utilizes a similar approach, where a TMOP-based quality estimator makes
refinement and derefinement decisions based on both (i) the function values and
(ii) the quality of the mesh.

%=====================================================================
%=====================================================================
\section{Methodology for \hra} \label{sec_method}

In this section we describe our new approach for \hra, which uses TMOP for \ra\
as discussed in Section \ref{subsec_ra}, and introduces a TMOP-based quality
estimator for \ha. The goal of this estimator is to determine the
refinement/derefinement decision that minimizes the TMOP objective function
$F$, i.e., the difference between the current ($A$) and target ($W$)
transformations.

\begin{figure}[t!]
\centerline{
\includegraphics[width=0.90\textwidth]{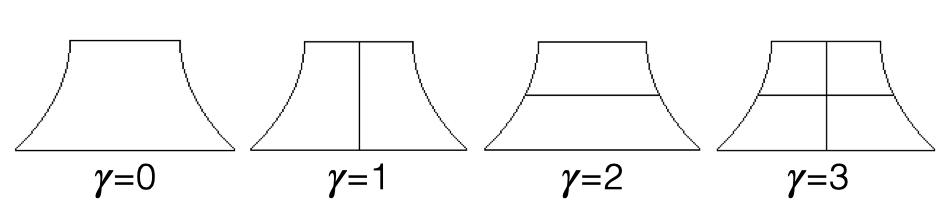} \hfil}
\caption{Example of the impact that different types of refinement can have on
a high-order curvilinear element in terms of size, aspect-ratio and skewness.}
\label{fig_ho_skew}
\end{figure}

Figure \ref{fig_amr_refs} shows how different types of refinement impact the
size and aspect-ratio of quadrilaterals and hexahedra. Isotropic refinement of
an element reduces its size and anisotropic refinement changes both the
aspect-ratio and the size of an element. Thus, \ha\ can be used with a
TMOP-based quality estimator to satisfy the size and aspect-ratio
\emph{targets} in the domain.  Note that for curved meshes, $h-$refinement of
an element could produce child elements that have different skewness.  Example
of such configurations can be observed in  Figure \ref{fig_ho_skew}.  As
described in the following paragraphs, the TMOP-based estimators naturally take
into account both the size and the shape changes.

%% TK: Maybe plot here aspect and skewness before and after refinement to emphasize the difference?

%==========================
\subsection{Refinement decision} \label{sec_refine}
Since refining one element does not affect the $A$ matrix of its neighbors, the
refinement decision can be performed independently in each element.  We define
$\Gamma$ as the set of all possible refinements that have to be considered for
\ha\ in a given element. If a \sz\ metric is used for the refinement decision,
only isotropic refinement is considered. If a \sh\ metric is used, only
anisotropic refinements are considered. For a \shsz\ metric, all refinements
are considered.  Once $\Gamma$ is constructed, our goal is to find the
refinement $\gamma \in \Gamma$ for each element that reduces its TMOP objective
function the most:
\begin{equation}
\label{eq_F_E2}
  \max_{\gamma} \Delta F^{\gamma}_E,\,\,\, \gamma \in \Gamma,
\end{equation}
where
\begin{eqnarray}
\label{eq_F_E}
  \Delta F^{\gamma}_E = F^{\gamma=0}_E - \frac{F^{\gamma}_{E}}{N_c}.
\end{eqnarray}
In \eqref{eq_F_E}, the first term on the right hand side denotes the value of
the objective function associated with the \emph{parent} element ($\gamma=0$),
and the second term denotes the mean of the objective function associated with
the corresponding $N_c$ \emph{children} elements obtained using refinement type
$\gamma$ on element $E$:
\begin{eqnarray}
\label{eq_F_E0}
  F^{\gamma}_E = \sum_{i=1}^{N_c} F_{E_c}.
\end{eqnarray}
If none of the refinements reduce the TMOP objective function for an element,
i.e.  $\Delta F^{\gamma}_E < 0 \,\,\forall\,\, \gamma$, then that element is
not refined.

%==========================

\subsection{Derefinement decision} \label{sec_derefine}

If the mesh contains children elements, i.e., elements that have been obtained
via refinement at a previous iteration of \hra, we consider them for
derefinement if restoring their parent element reduces the TMOP objective
function. The derefinement capability is crucial for time-dependent problems
where the shape and size targets can change with time as the solution evolves
or if the original mesh has more resolution than required to begin with.  Using
$E_p$ to denote a parent element that has been refined to spawn $N_c$ children
elements at a previous iteration of \hra, the parent element is derefined if
$\Delta F_{E_p} > 0$, where
\begin{equation}
\label{eq_F_E3}
  \Delta F_{E_p} = \sum_{i=1}^{N_c} \frac{F_{E_c}}{N_c} - F_{E_p}.
\end{equation}

As we can see, \eqref{eq_F_E} and \eqref{eq_F_E3} are compliments of each
other.  In \eqref{eq_F_E}, we compare how the refinement $\gamma$ changes the
TMOP functional between the parent and its children element, and in
\eqref{eq_F_E3}, we see how removing a refinement done using \eqref{eq_F_E} in
the past changes the TMOP functional.

%==========================

\subsection{Algorithm for \hra}
The TMOP-based \hra\ method is summarized in Algorithm \ref{algo_hra}.  Here,
the input consists of the original mesh positions ($\bx_0$), the mesh quality
metric ($\mu$), the tolerance ($\epsilon$) that depends on the norm of the
gradient of the TMOP objective function at the current nodal coordinates, and
the number of $h-$refinements ($N_h$) to be done after each $r-$refinement.
The output is the optimized mesh positions $\bx_s$.  The algorithm performs
consecutive \ra\ and \ha\ steps until the \ha\ step does not result in any
changes in the mesh output by \ra. Thus, using $N_R$ and $N_D$ to denote the
total number of elements that are refined and derefined, respectively, by the
\ha\ component (Step 6), the \hra\ procedure stops when $N_R = N_D = 0$.

\begin{algorithm}[h!]
\label{algo_hra}
\SetAlgoLined
\KwIn{$\bx_0$, $\mu$, $\epsilon$, $N_h$.}
\KwOut{$\bx_s$ (initialized to $\bx_0$).}
Construct $W_i$ for each integration-point $i$ using target-construction. ~\cite{knupp2019target}\\
\While{$N_R\ne 0 \,\,{\tt or }\,\, N_D \ne 0$}{
    $\text{\ra}$:\\
    $\bx_s \rightarrow \underset{\bx}{\argmin} \,\,\,
  \sum_{E \in \mathcal{M}} \sum_{\bx_q \in E_t}
                           w_q\,\det(W(\bar{\bx}_q))\, \cmb(\bx_q) \mu(T(\bx_q))$. ~\cite{TMOP2020}\\
    $\text{\ha}$:\\
    \For{$i \in 1 \dots N_h$} {
          $\forall\,\,\, E_p\,\,\, \in\,\,\, \mathcal{M}$, determine $\Delta F_{E_p}$. ~\eqref{eq_F_E3}\\
          Derefine element $E_p$ if $\Delta F_{E_p} > 0$. \\
          $\forall\,\,\, E\,\,\, \in\,\,\, \mathcal{M}$, determine $\Delta F^{\gamma}_E\,\,\, \forall\,\,\, \gamma\,\,\, \in\,\,\, \Gamma$. ~\eqref{eq_F_E2} \\
          Refine element $E$ based on $\max_{\gamma} \Delta F^{\gamma}_E$. ~\cite{cerveny2019nonconforming}
    }
}
\caption{\hra}
\end{algorithm}

In Algorithm \ref{algo_hra}, determining $\max_{\gamma} \Delta F^{\gamma}_E$
requires each element to be refined using each refinement type in $\Gamma$
followed by the targets $W$ and $A$ to be reconstructed on the quadrature
points of the children element. This remap of targets from parent to children
element is the main source of increase in computational cost for \hra\ in
comparison to \ra.

The flexibility of the \hra\ method can be further increased by utilizing
different mesh quality metrics for \ra\ and \ha\ in Steps 3 and 5,
respectively.  This allows us to use different quality metrics to optimize
different aspects of the mesh using \ha\ ($\mu^h$) and \ra\ ($\mu^r$).  For
example, using $\mu^h_2$ (\emph{shape} metric) and $\mu^r_7$
(\emph{shape}$+$\emph{size} metric) allows us to optimize the aspect-ratio with
\ha, and skew, aspect-ratio and size with \ra.

%==========================

\subsection{Target reconstruction during mesh optimization}
For mesh adaptivity, the targets are either defined through analytical
functions or through discrete functions that depend on the discrete solution
defined on the original mesh. For analytical function based adaptivity, it is
straightforward to construct targets at each degree of freedom as the mesh
moves during \ra\ or as new elements are added during \ha. For discrete
adaptivity however, the discrete functions must be mapped from the original
mesh ($\bx_0$) to the active mesh ($\bx$) so that the targets can be
reconstructed to determine $\mu(T)$.

\begin{figure}[h!]
\centerline{
  \includegraphics[width=0.40\textwidth]{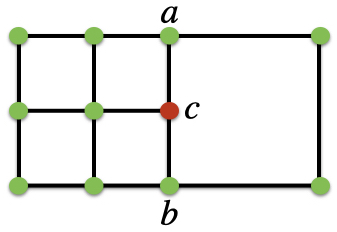} \hfil}
\caption{Isotropic refinement of an element in the mesh showing the
true degrees of freedom (colored green) and the constrained degree
of freedom (labeled ``c''). The conforming prolongation operator
defines the value at ``c'' as the average of the values in ``a'' and ``b''.}
\label{fig_amrdonor}
\end{figure}

During \ra, the discrete functions are remapped using an advection-based PDE,
as described in Section 4.2 of \cite{IMR2018}, or using high-order
interpolation between the meshes (see Section 2.3 of
\cite{mittal2019nonconforming}).  For \ha, mesh refinement introduces elements
that are a subset of existing elements.  Since every discrete function defined
on the mesh is represented as a polynomial on each element, it is
straightforward to interpolate the desired discrete function from a given
parent element to its children. The key challenge in mapping the discrete
function from the original mesh to the refined mesh is maintaining the desired
continuity in the solution at the nonconforming interfaces of the element
boundaries.  This continuity of the original discrete solution is maintained
using a conforming prolongation operator, as described in
\cite{cerveny2019nonconforming}, that extends any function described on the
(unconstrained) \emph{true} degrees of freedom to the (constrained) hanging
nodes resulting from refinement operations on the original mesh. Figure
\ref{fig_amrdonor} shows an example of a nonconforming mesh where isotropic
refinement on an element introduces a constrained degree of freedom (labeled
``c'') that depends on the true degrees of freedom (labeled ``a'' and ``b'').

%==========================

\subsection{Mesh conformity during \ra}
An important consideration for the \ra\ method described in this work is
ensuring that the optimal nodal positions determined in Step 3 of Algorithm
\ref{algo_hra} maintain continuity between the nonconforming interfaces of the
mesh.

In our variational formulation \eqref{eq_F_full}, the TMOP objective function
is computed as an element-by-element weighted sum of the quality metric that
measures the deviation between the current and target geometrical parameters of
the mesh. Consequently, the derivatives of the TMOP objective function required
for \ra\ are also computed on an element-by-element basis with contribution
from each nodal finite element degree of freedom.  To maintain continuity
between the nonconforming interfaces of the mesh, we project the gradients from
the local nodes to the true degrees of freedom using the transpose of the
projection operator. Thus, the optimization solver (e.g., Newton's method)
operates on the true degrees of freedom to compute optimum nodal positions,
which are used to update the positions of the constrained nodes using the
prolongation operator. This approach ensures that the nonconforming interfaces
of the mesh stay consistent during nodal movement with \ra.

\subsection{Alternative approaches for TMOP estimator during \ha}
An alternative to determining how a refinement type changes the TMOP objective
function, without explicitly refining each element, is to decompose $A$ into
its four geometrical components (e.g., \eqref{eq_W}), and then scale the
geometrical components based on the refinement type being considered.  For
example, to estimate how $\gamma=1$ would change the TMOP objective function in
a 2D quadrilateral, we take $A$ associated with the quadrature points of the
sample point and scale its volume component ($\zeta$) by a factor of 2 to
capture the size-reduction and scale the diagonal entries of $D$ to indicate an
increase in aspect-ratio in $x$ by a factor of 2.  Using this approach, the
TMOP estimator becomes:
\begin{equation}
% \label{eq_F_E3}
  \Delta F^{\gamma}_E = F^{\gamma=0}_E - \tilde{F}^{\gamma}_E,
\end{equation}
where $\tilde{F}^{\gamma}_E$ represents the functional associated with the
parent element after the geometrical components of $A$ have been scaled based
on $\gamma$.  This approach is computationally cheaper than the approach in
Algorithm \ref{algo_hra} because it does not require each element to be
explicitly refined and the targets $W$ to be mapped from the parent element to
its children.  Numerical experiments show that this approach is as effective as
the approach in Algorithm \ref{algo_hra} when a low-order mesh is optimized.
This behavior is expected because rescaling only the size and aspect-ratio
components in \eqref{eq_W} to determine the impact of refinement type ignores
its impact on the skewness of an element. Since our goal is the optimization of
curvilinear high-order meshes, we forego this low-order approach due to
higher-accuracy of the approach in Algorithm \ref{algo_hra}.

Considering the modular structure of the \hra\ framework, the TMOP-based
estimator for $h-$refinement can also be replaced or used in conjunction with
other (error) estimators if desired.  All numerical experiments presented in
the current work use the TMOP-based estimator.

%=====================================================================
%=====================================================================

\section{Numerical experiments} \label{sec_apps}
In this section, we present several numerical experiments to illustrate the
effectiveness of the \hra\ method and compare it with standalone \ra\ and \ha\
in satisfying the geometric parameters' targets.  The new algorithms were
implemented in the open-source MFEM finite element library
\cite{mfem,mfem-web}.

%==========================

\subsection{2D benchmark using the Poisson equation} \label{subsec_pois}
\begin{figure}[tb!]
\begin{center}
$\begin{array}{cc}
\includegraphics[height=0.25\textwidth]{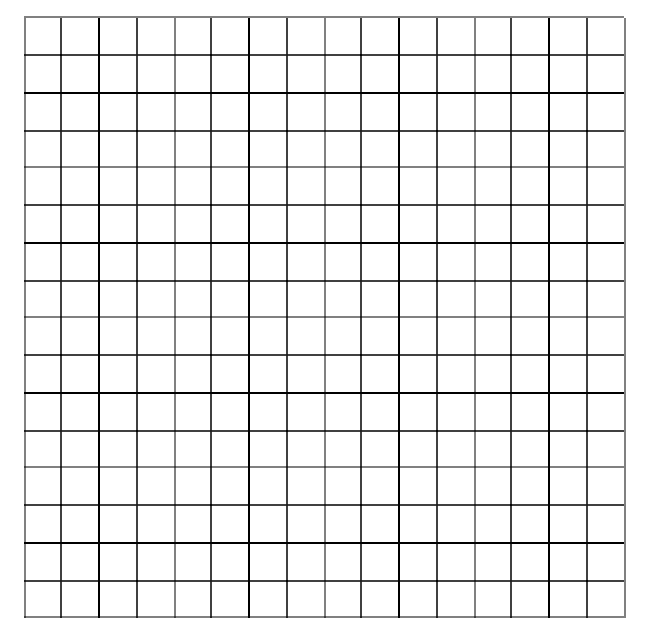} &
\includegraphics[height=0.25\textwidth]{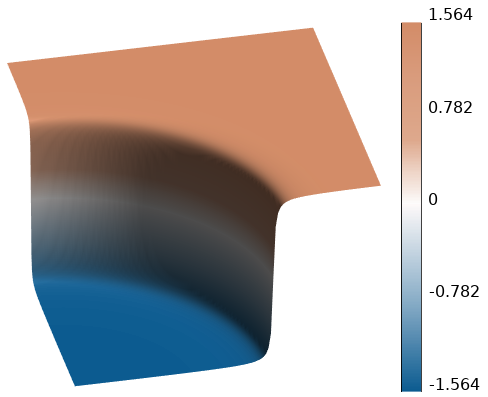} \\
\end{array}$
\end{center}
\caption{(left) Original $16\times 16$ mesh and (right) the exact solution to the Poisson problem
 considered here \eqref{eq_pois_u}.}
\label{fig_pois}
\end{figure}

To quantify the improvement in the accuracy of the solution due to \hra\ over
\ra\ and \ha, we consider the 2D benchmark from \cite{cerveny2019nonconforming}
where the Poisson equation with a known exact solution is solved in $\Omega=[0,
1]^2$,
\begin{equation}
\label{eq_pois}
\nabla^2 u = f.
\end{equation}
Here $f$ is chosen such that
\begin{equation}
\label{eq_pois_u}
u = \arctan\bigg[\alpha\bigg(\sqrt{(x-x_c)^2 + (y-y_c)^2}\bigg)\bigg],
\end{equation}
which has a sharp circular wave front of radius $r$ center around $(x_c, y_c)$.
Figure \ref{fig_pois} shows the original second order ($Q_2$) mesh with the
exact solution for $r = 0.7$ and  $\alpha = 200$.  Below we illustrate the
effectiveness of different mesh-adaptivity techniques such as $h-$, $r-$, and
\hra\ in reducing the error in the solution without having to remesh the
domain.

Since we use a \shsz\ metric ($\mu_9$) to adapt the mesh, the target elements
must include information about both shape and size.  We use the exact solution
to construct the target transformation,
\begin{eqnarray}
\label{eq_pois_W}
W =
\sqrt{\zeta}
\begin{bmatrix}
1 & 0 \\
0 & 1
\end{bmatrix}
\begin{bmatrix}
1 & \cos\,\phi \\
0 & \sin\,\phi
\end{bmatrix}
\begin{bmatrix}
\frac{1}{\sqrt{\rho}} & 0 \\
0 & \sqrt{\rho}
\end{bmatrix},
\end{eqnarray}
where the size target ($\zeta$) depends on the magnitude of the gradient (i.e.,
$\zeta \propto ||\nabla u||$), target skewness ($\phi$) is set to be the same
as that for an ideal element, i.e., $\phi = \pi/2$, and aspect ratio target
($\rho$) is computed by the ratio of the gradient components (i.e., $\rho
\propto \nabla_x u/\nabla_y u$).  Thus,
\begin{eqnarray}
\label{eq_pois_W2}
W =
\begin{bmatrix}
\sqrt{\frac{\zeta}{\rho}} & 0 \\
0 & \sqrt{\zeta\rho}
\end{bmatrix}.
\end{eqnarray}

\begin{figure}[tbh!]
\begin{center}
$\begin{array}{ccc}
\includegraphics[height=0.25\textwidth]{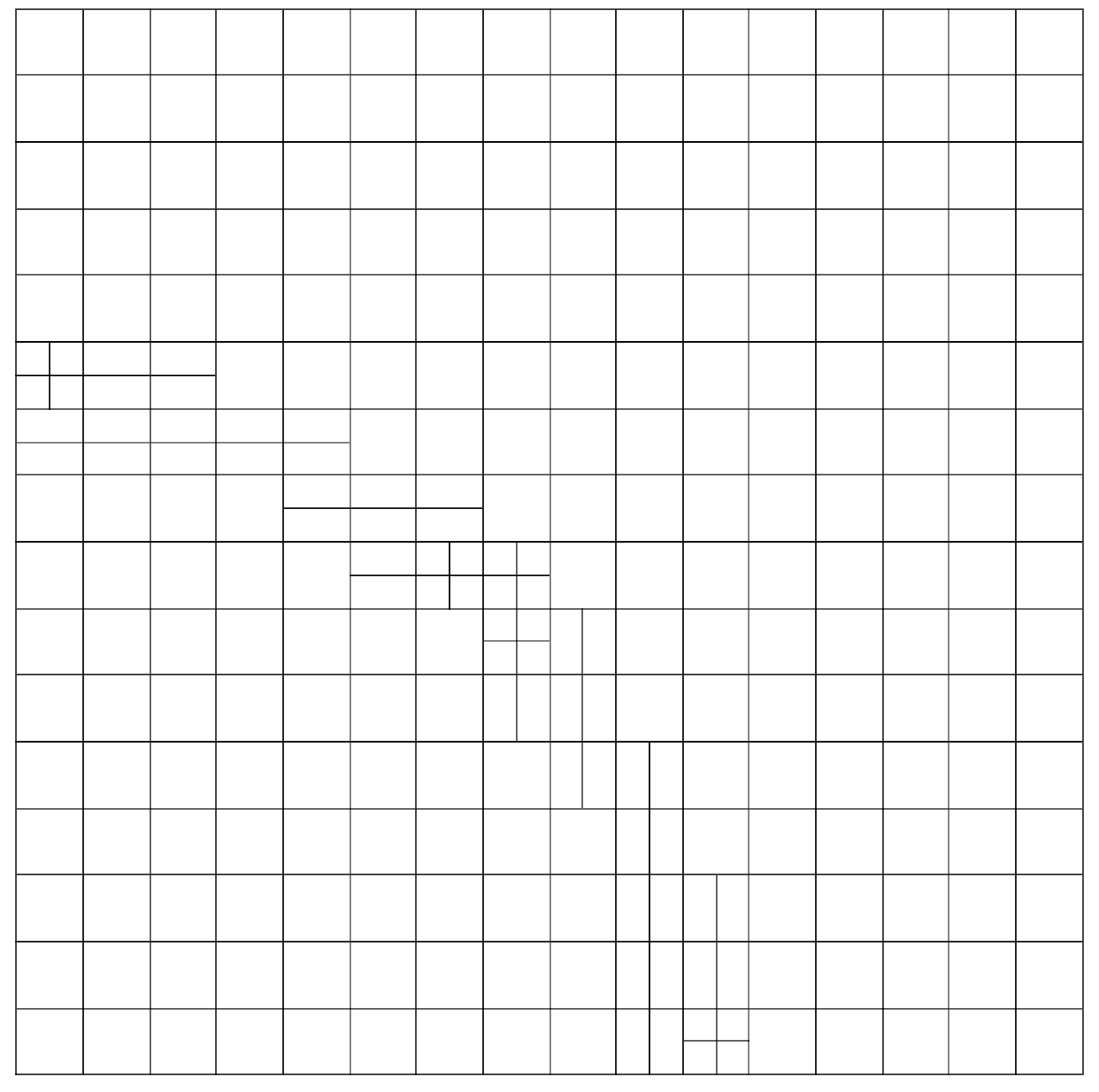} &
\includegraphics[height=0.25\textwidth]{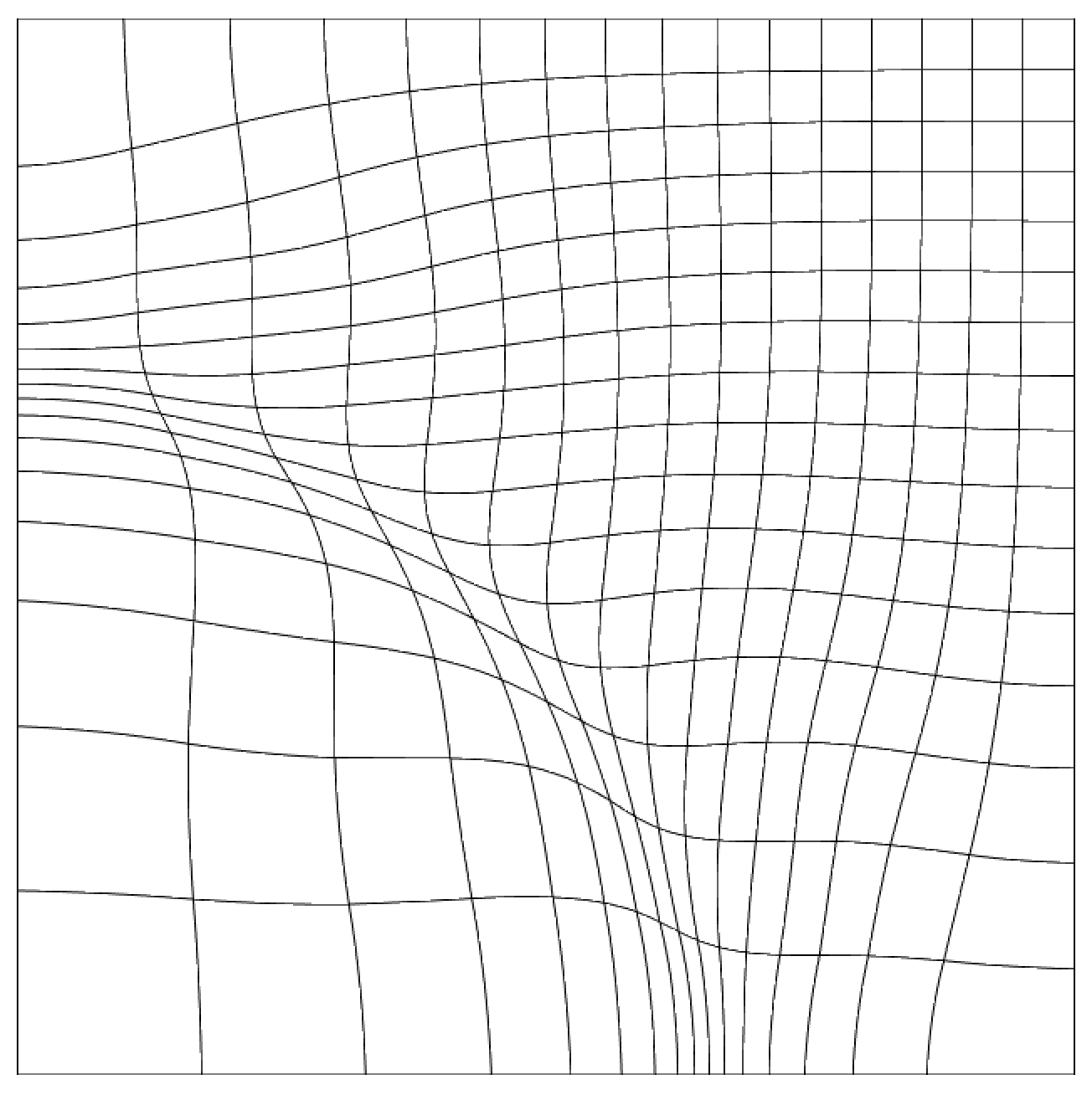} &
\includegraphics[height=0.25\textwidth]{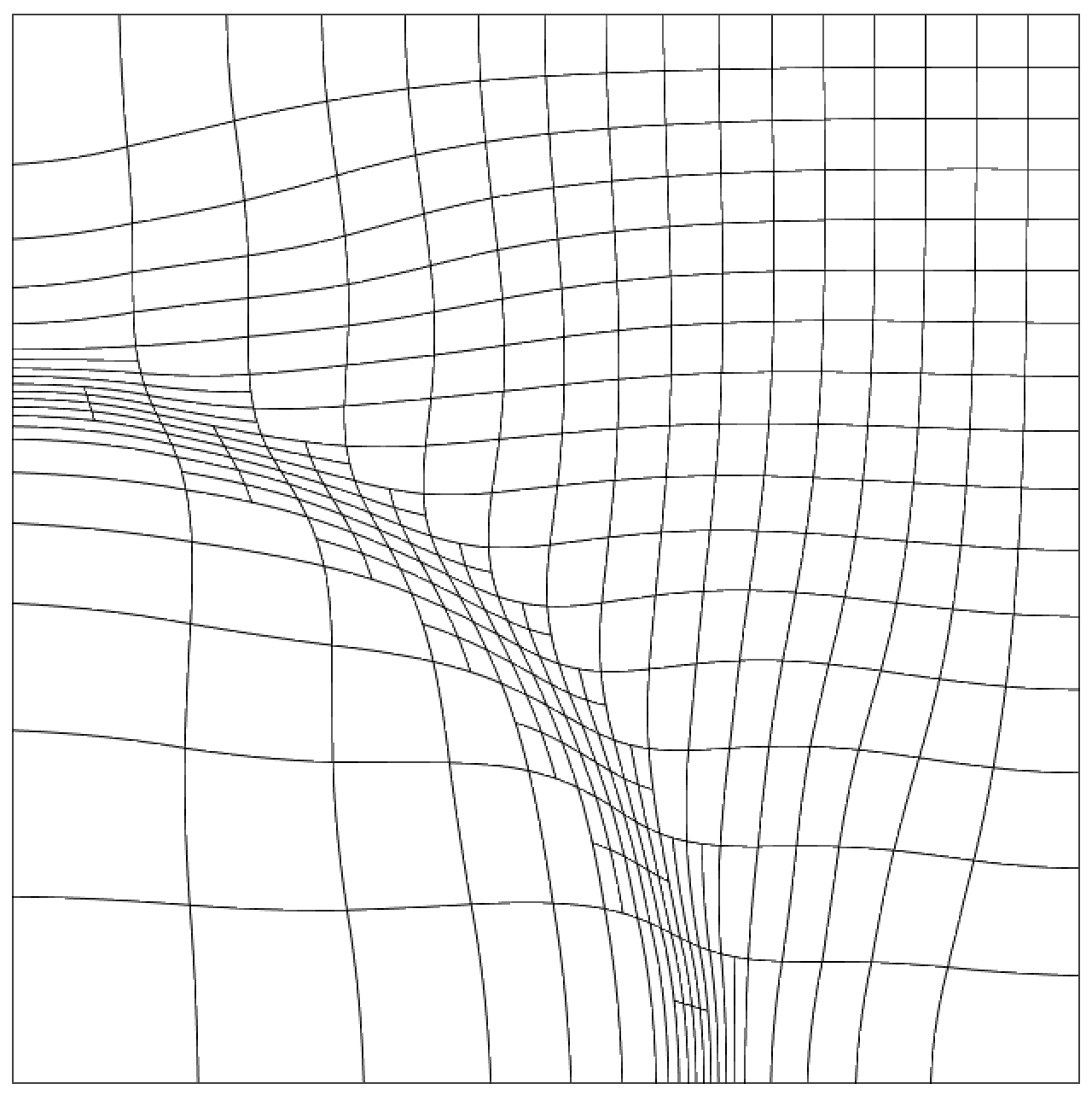} \\
\end{array}$
\end{center}
\vspace{-7mm}
\caption{Optimized meshes obtained using (left) \ha, (center) \ra, and
(right) \hra.}
\label{fig_pois_optimized}
\end{figure}

Figure \ref{fig_pois_optimized} shows optimized meshes obtained by $h$-, $r$-,
and \hra.  The initial mesh for these simulations is the one from Fig.
\ref{fig_pois}.  As evident, each of these three methods increases the mesh
resolution in the region with sharp solution gradients. While \ra\ moves nodal
positions to increase the resolution, \ha\ uses anisotropic refinement near the
domain boundaries (to satisfy aspect-ratio targets) and isotropic refinement
away from the boundaries (to satisfy size targets).

\begin{figure}[b!]
\begin{center}
$\begin{array}{c}
\includegraphics[height=0.5\textwidth]{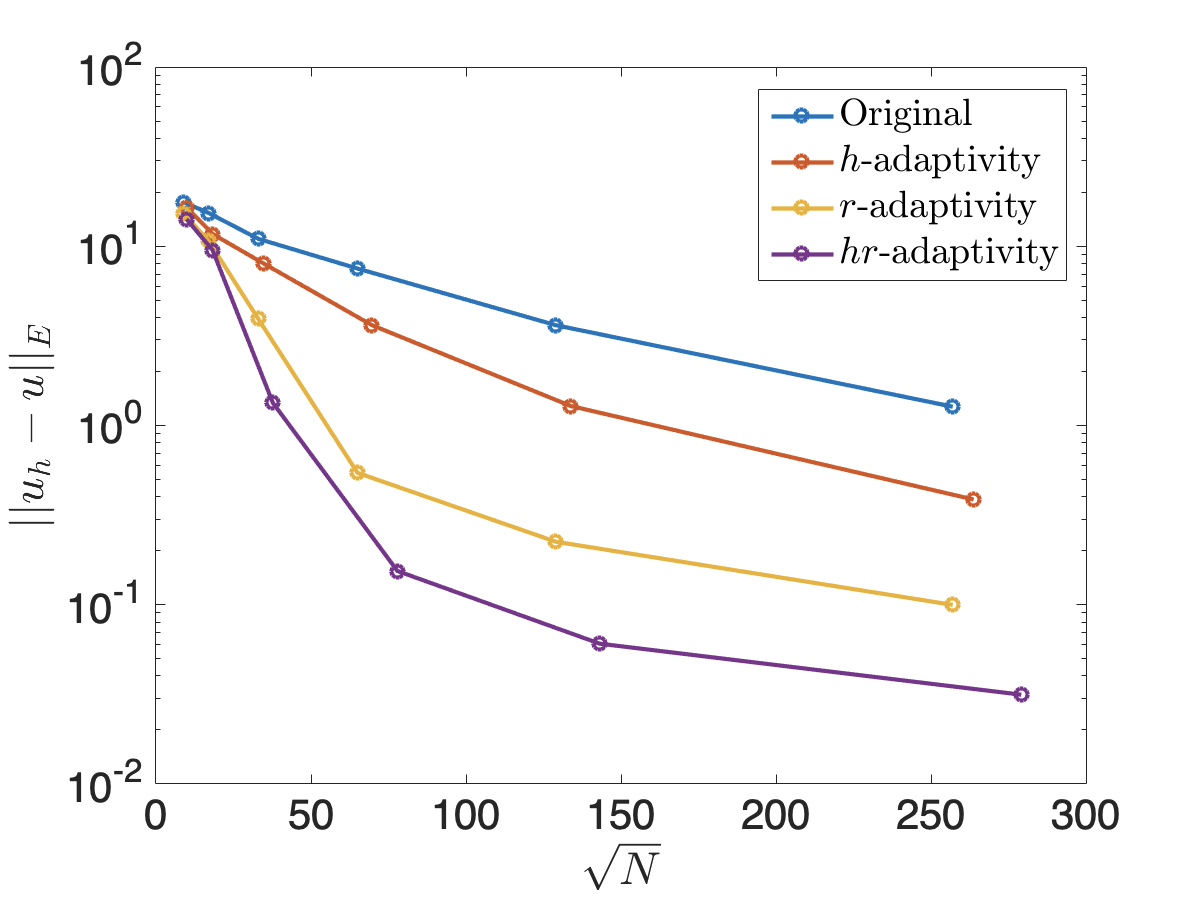} \\
\end{array}$
\end{center}
\caption{Error for the discrete solution ($u_h$) with respect to the exact
solution ($u$) in the energy norm for a sequence of uniform meshes (the
original mesh) and meshes obtained
using \ha, \ra, and \hra.}
\label{fig_pois_error}
\end{figure}

To obtain a better quantitative comparison, we start with a $4\times 4$ mesh
and increase the resolution by a factor of 2 in each direction, while adapting
the mesh using \ha, \ra, and \hra\ for each case.  For \ha\ we use the same
criterion for refinement that is described in Algorithm \ref{algo_hra}.
Additionally, for \ha\ and \hra, we consider only 1 iteration of each kind of
mesh refinement to illustrate the effectiveness of these methods while keeping
the added computational cost to a minimum.  Once the mesh has been optimized,
the Poisson problem is solved and the error is measured for the discrete
solution ($u_h$) with respect to the exact solution ($u$) in the energy norm
\cite{cerveny2019nonconforming}.  Figure \ref{fig_pois_error} compares the
error for a sequence of uniform meshes (the ``original'' meshes) with meshes
obtained using \ha, \ra, and \hra. In Fig. \ref{fig_pois_error}, $N$ represents
the total number of degrees of freedom in the mesh. As evident, \ra\ can help
significantly reduce the error in the solution by moving the nodal positions of
the mesh and even a single iteration of $h-$refinement with \hra\ is sufficient
for an even more significant improvement. This numerical experiment shows that
we require about 66\% fewer degrees of freedom with \hra\ as compared to \ra\
for a given accuracy in the solution. Note that one can do multiple iterations
of \ha\ and \ra\ during the mesh-optimization process, but we do a single
iteration here to illustrate the effectiveness of these techniques in improving
the accuracy of the solution while keeping the computation cost of the
mesh-optimization process to a minimum.

%==========================

\subsection{Analytical adaptivity} \label{subsec_analytic}
The next example that we consider is adaptivity based on analytical functions
to demonstrate how \hra\ can be used for purely \sz-adaptivity and for
\shsz-adaptivity. The \ra\ component in each case is based on a \shsz\ metric,
and we setup the targets to mimic a problem with a shock wave propagating
externally from the center of the domain.

For \hra\ with a \sz\ metric, the target $W$ is defined in the domain $\Omega
\in [0,1]^2$ such that the target skewness is the same as that for an ideal
element, i.e., $\phi = \pi/2$, target aspect-ratio is unity, i.e. $\rho = 1$,
and target size ($\zeta$) depends on an analytic function $\eta \in [0,1]$ that
is a function of physical-space coordinates ($\bx$). Thus,
\begin{eqnarray}
\label{eq_W_ex1}
W =
\begin{bmatrix}
\sqrt{\zeta} & 0 \\
0 & \sqrt{\zeta}
\end{bmatrix},
\end{eqnarray}
where
\begin{equation}
\label{eq_W_ex1_eta}
\eta = \tanh(\beta(r-0.2)) - \tanh(\beta(r-0.3)).
\end{equation}
In \eqref{eq_W_ex1_eta}, $\beta$ determines the sharpness in the gradient of
the solution, and $r$ is the distance from the center of the domain.  Using
$\eta$, the size targets are defined as:
\begin{equation}
\label{eq_W_ex1_delta}
\zeta = \eta \delta + (1 - \eta) \psi \delta,
\end{equation}
where $\delta$ is the target size of the degrees of freedom in the region with
$\eta=1$ and $\psi$ is the ratio of the biggest to smallest elements in the
domain. For the example considered here, we set $\beta=30$, $\delta = 0.001$,
and $\psi=10$.

\begin{figure}[t!]
\begin{center}
$\begin{array}{cccc}
\includegraphics[height=25mm]{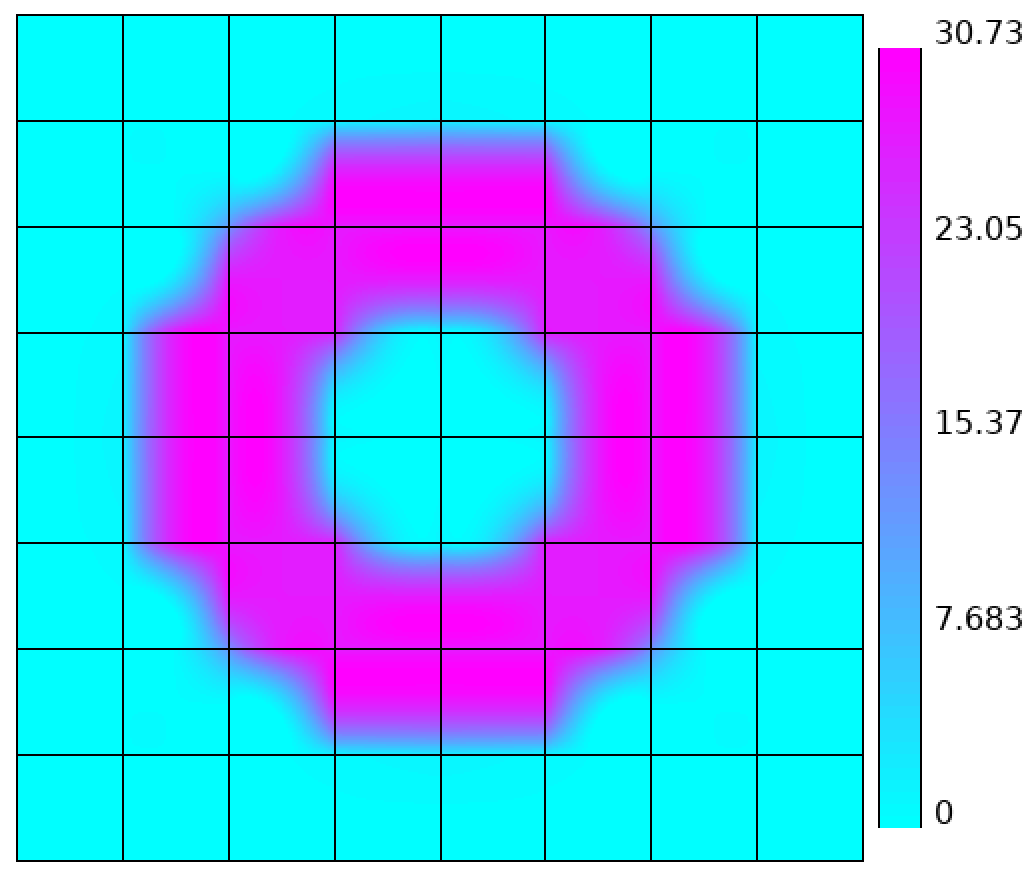} &
\includegraphics[height=25mm]{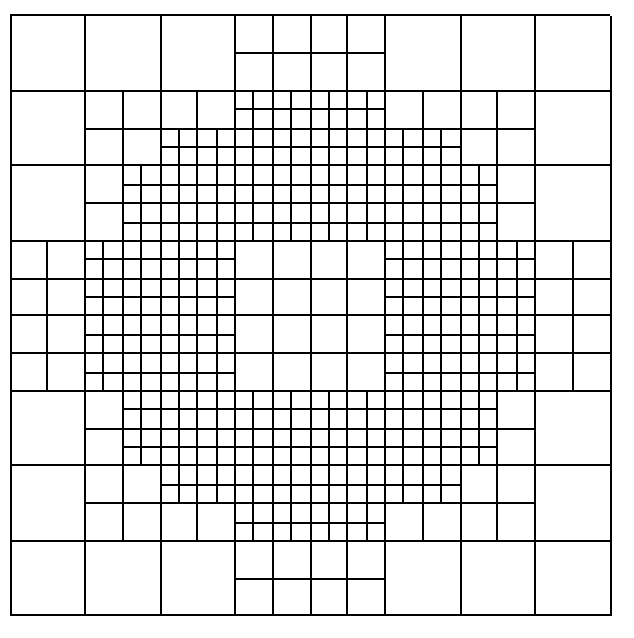} &
\includegraphics[height=25mm]{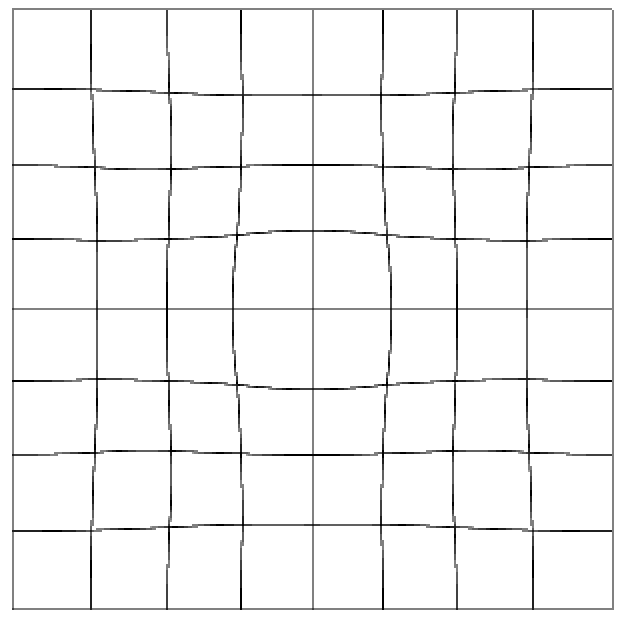} &
\includegraphics[height=25mm]{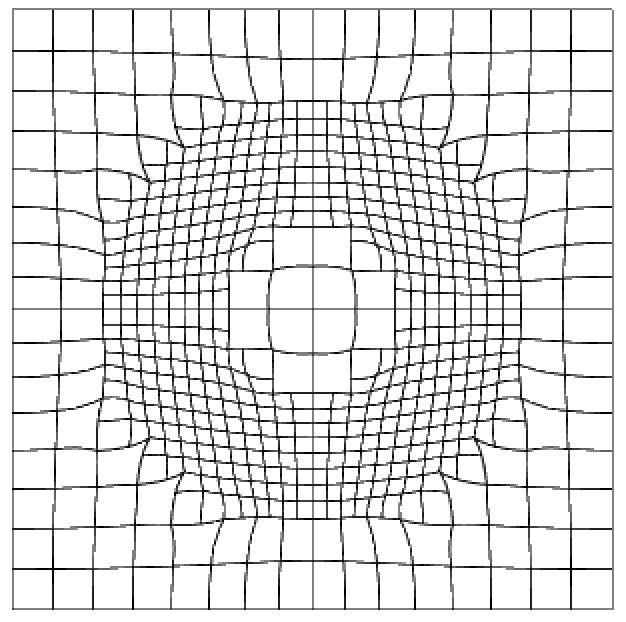} \\
\multicolumn{4}{c}{\textrm{(i) $8\times8$ mesh}} \\
\includegraphics[height=25mm]{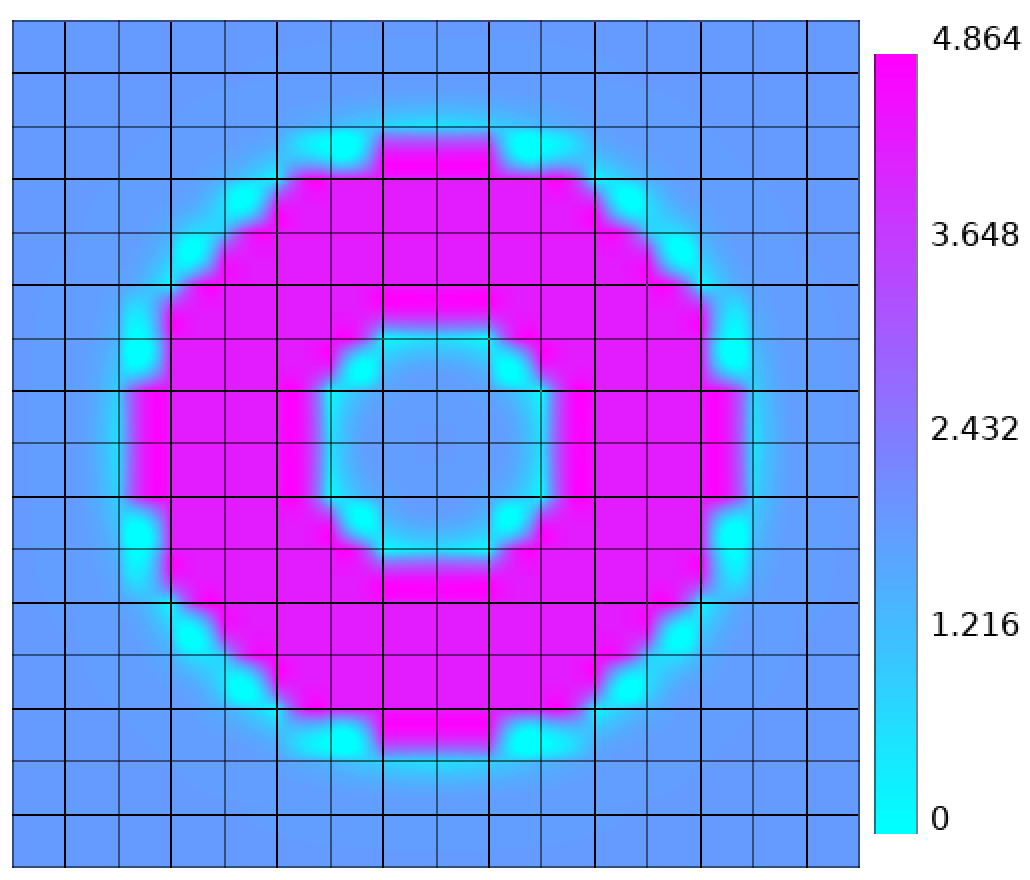} &
\includegraphics[height=25mm]{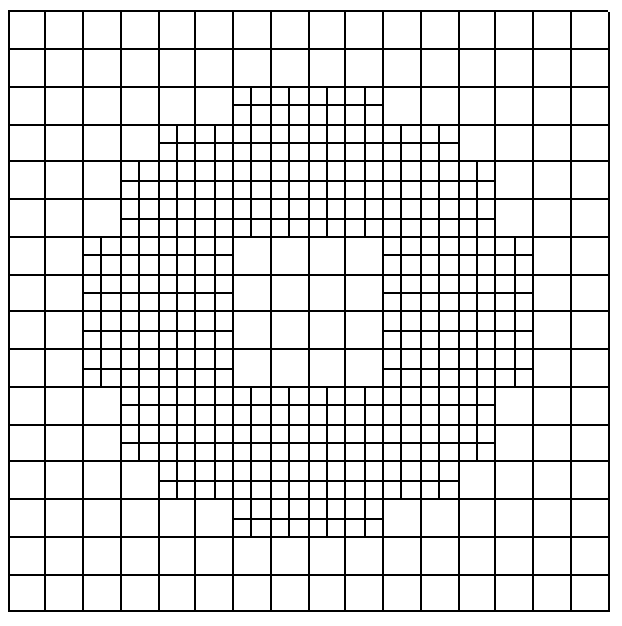} &
\includegraphics[height=25mm]{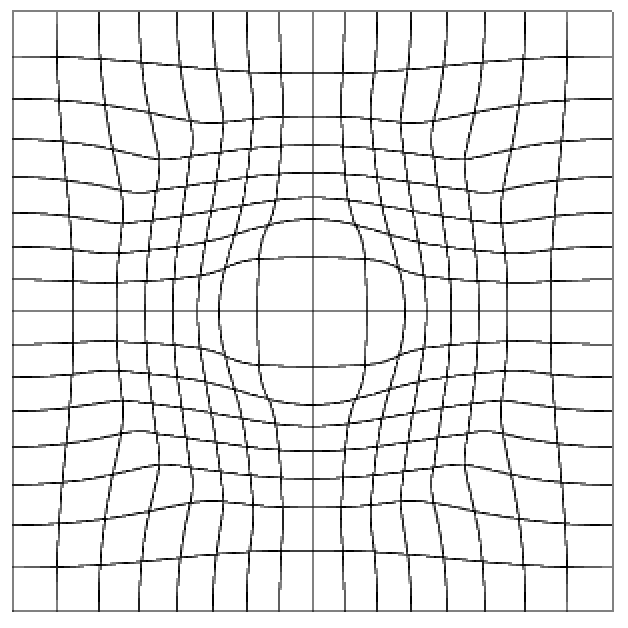} &
\includegraphics[height=25mm]{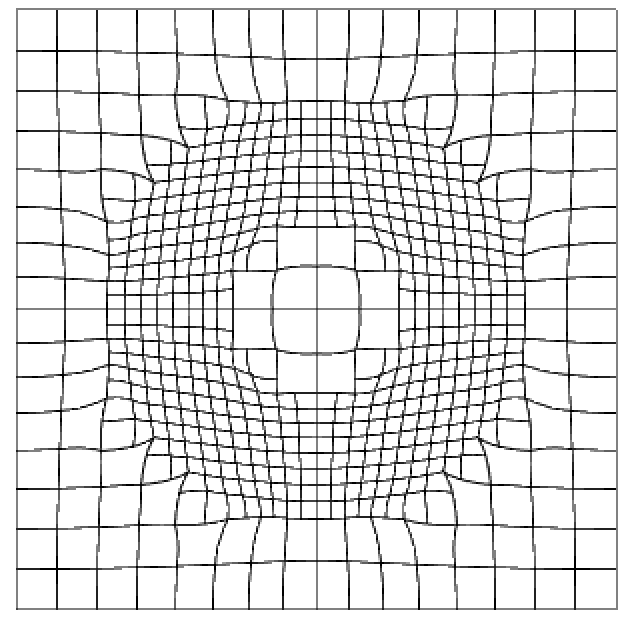} \\
\multicolumn{4}{c}{\textrm{(ii) $16\times16$ mesh}}
\end{array}$
\end{center}
\vspace{-7mm}
\caption{Analytic \sz\ target for optimizing (i) $8\times8$ and (ii)
$16\times16$ meshes with equally sized elements. In each case, we show (left to
right) (a) the original mesh with metric values $\mu^r(T)$ compared with meshes
obtained using (b) \ha, (c) \ra, and (d) \hra.}
\label{fig_app_ex1}
\end{figure}

Figure \ref{fig_app_ex1} shows examples of the difference in the optimized
meshes obtained using \hra\ ($\mu^h_{55}$ and $\mu^r_7$) versus \ra\
($\mu^r_7$) and \ha-only ($\mu^h_{55}$) for two different initial mesh
resolutions.  For each case, the left panel shows the quality metric
$\mu^r_7(T)$ evaluated at the degrees of freedom of the original mesh, which
indicates the difference between the current and the target geometric
parameters \eqref{eq_W_ex1}.  For the mesh with $N_E=64$ equally sized
elements, \ha\ reduces the TMOP objective function $F(\bx)$ by 40.36\% while
increasing the element count to $N_E=484$, \ra\ does not reduce $F(\bx)$, and
\hra\ reduces $F(\bx)$ by 69.2\% while increasing the element count to
$N_E=616$.  For the $16\times16$ mesh, \ha\ reduces $F(\bx)$ by 21.9\% with
$N_E=544$, \ra\ reduces $F(\bx)$ by 51.8\%, and \hra\ reduces $F(\bx)$ by
67.3\% with $N_E=616$.

\begin{figure}[t!]
\begin{center}
$\begin{array}{cccc}
\includegraphics[height=25mm]{figures/app1_coarse_grid+metric} &
\includegraphics[height=25mm]{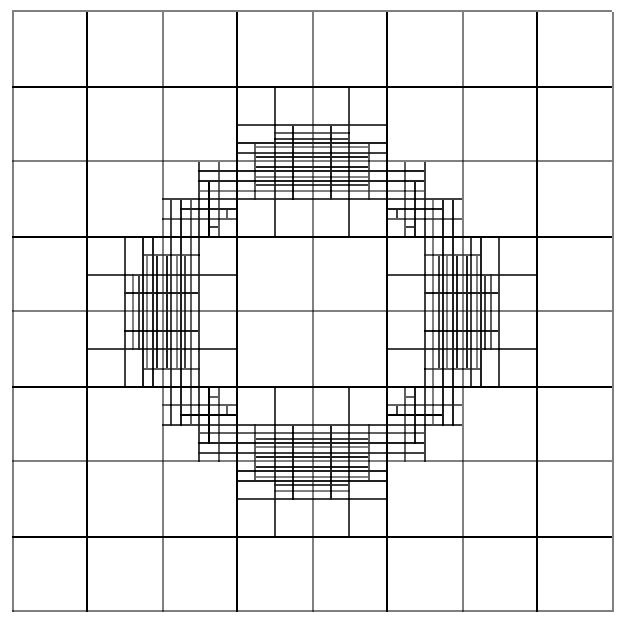} &
\includegraphics[height=25mm]{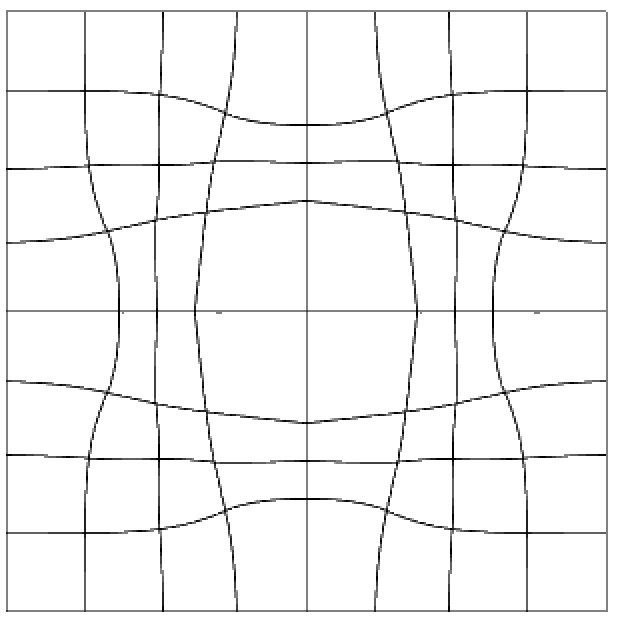} &
\includegraphics[height=25mm]{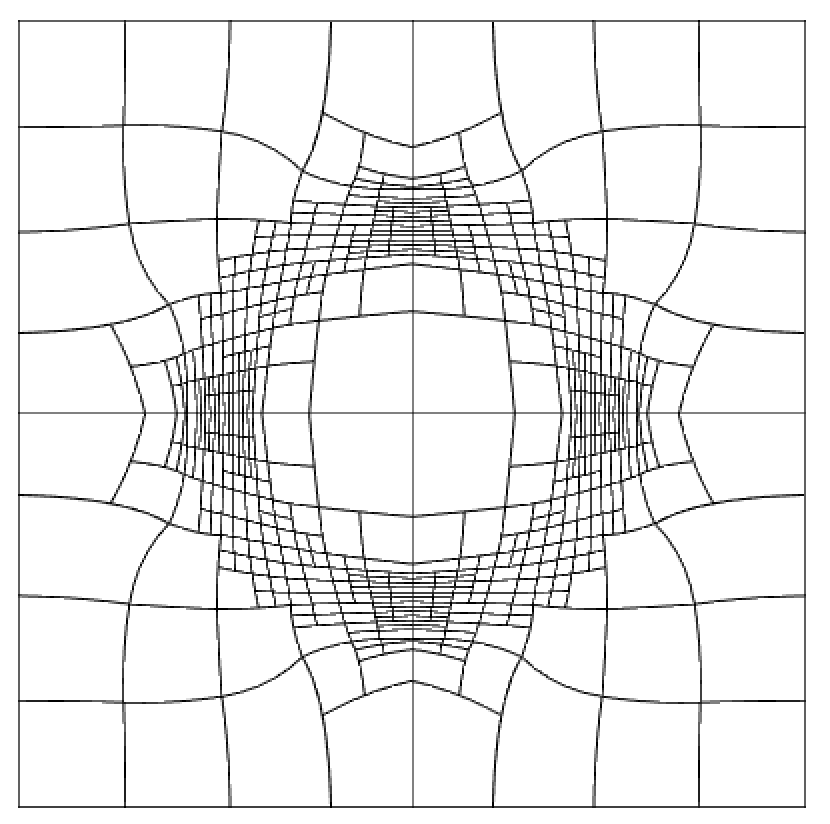} \\
\multicolumn{4}{c}{\textrm{(i) $8\times8$ mesh}} \\
\includegraphics[height=25mm]{figures/app1_fine_grid+metric} &
\includegraphics[height=25mm]{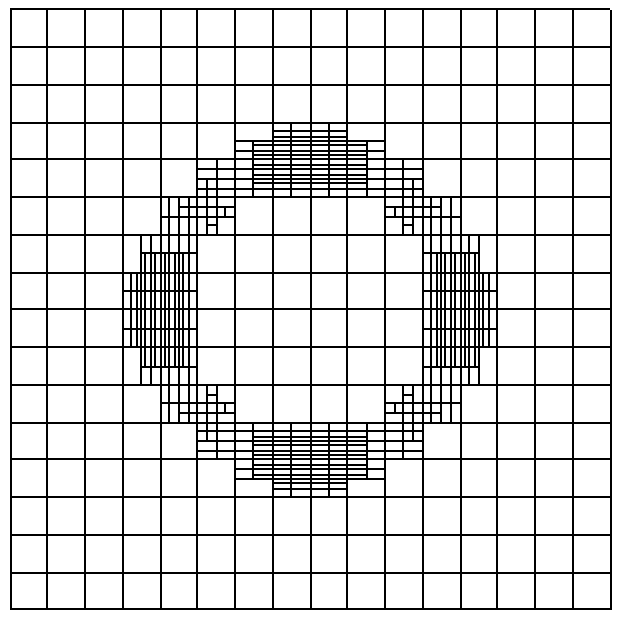} &
\includegraphics[height=25mm]{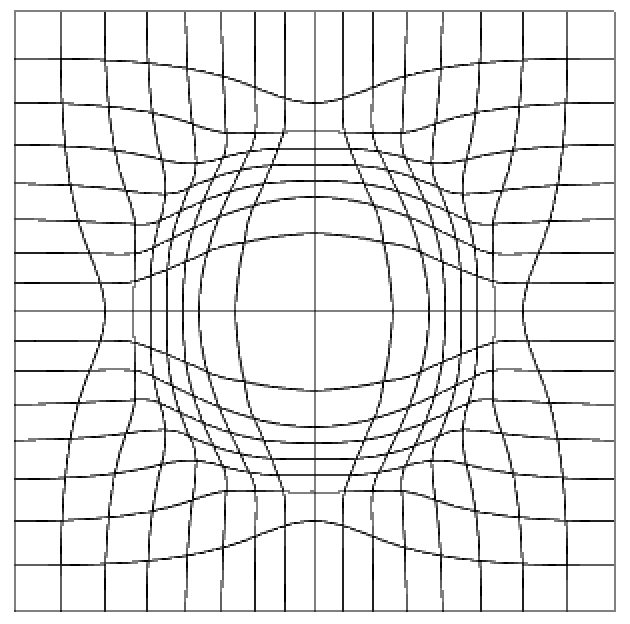} &
\includegraphics[height=25mm]{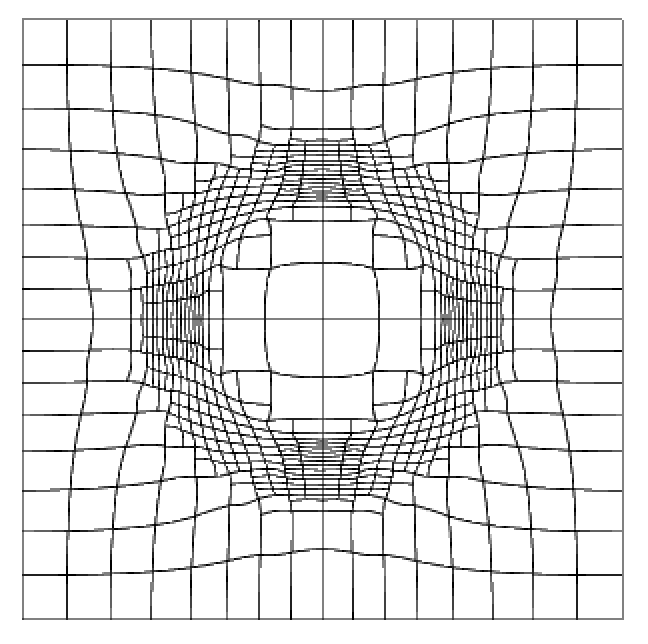} \\
\multicolumn{4}{c}{\textrm{(ii) $16\times16$ mesh}}
\end{array}$
\end{center}
\vspace{-7mm}
\caption{Analytic \shsz\ target for optimizing (i) $8\times8$ and
(ii) $16\times16$ meshes with equally sized elements. In each case, we show (left to right)
(a) the original mesh with metric values $\mu^r(T)$
compared with meshes obtained using (b) \ha, (c) \ra, and (d) \hra.}
\label{fig_app_ex2}
\end{figure}

For \hra\ with a \shsz\ metric ($\mu^h_{7}$), we modify the target ($W$) such
that the aspect-ratio depends on the distance from the center ($r$) and the
angle around the center of the domain, i.e., $\rho = f(r, \theta)$. Figure
\ref{fig_app_ex2} shows that \hra\ is able to optimize the mesh and satisfy the
\shsz\ targets much better as compared to \ra. For the $8\times8$ mesh, \ha\
reduces $F(\bx)$ by 92.7\% with a final element count of $N_E=572$, \ra\
reduces $F(\bx)$ by 6.1\%, and \hra\ reduces $F(\bx)$ by 93.5\% with $N_E=704$.
For the $16\times16$ mesh, $F(\bx)$ is reduced 65.6\% by \ha, 60.6\% by \ra,
and 85.3\% by \hra, with $N_E=600$ for \ha\ and $N_E=716$ for \hra.

%==========================

\subsection{Mesh adaptivity with simplices} \label{subsec_tri}
The \hra\ method described in this paper extends to simplices, where the
refinements are always isotropic.  So far, we have looked at examples with
structured quadrilateral meshes. Due to the ease with which more general
domains can be meshed with simplices, however, triangles/tetrahedra are often
preferred over quadrilaterals/hexahedra in complicated domains.

\begin{figure}[b!]
\begin{center}
$\begin{array}{cccc}
\includegraphics[height=35mm]{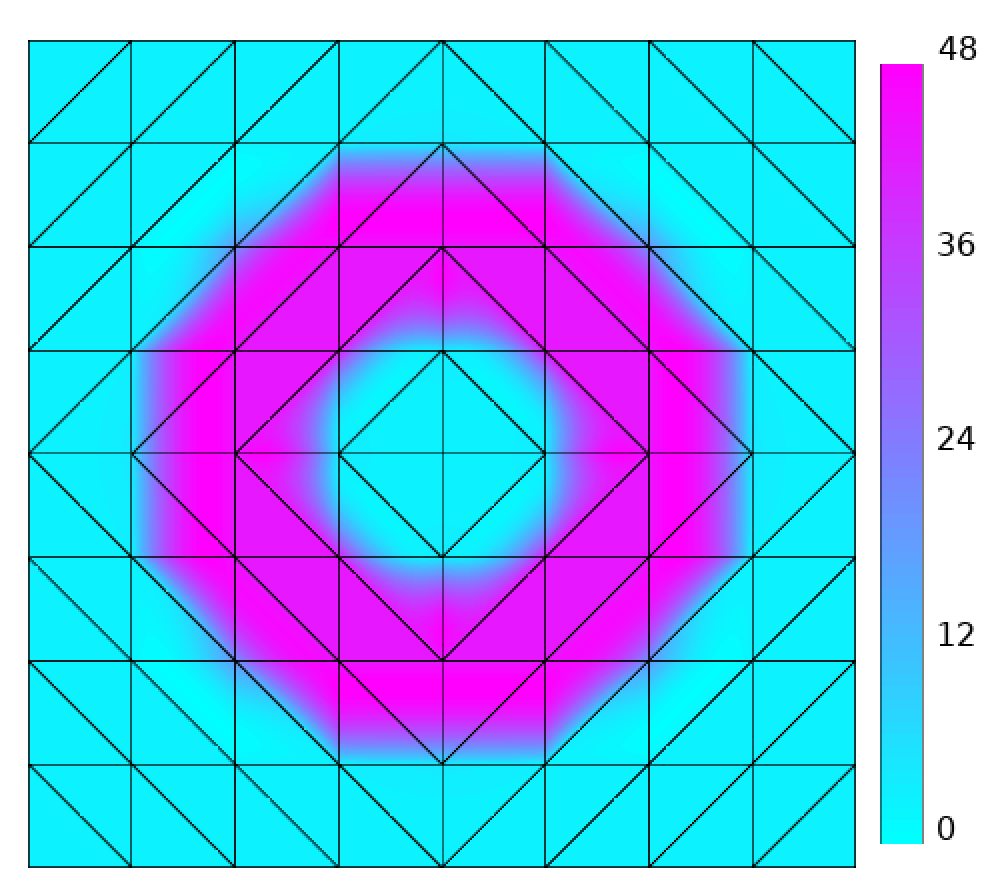} &
\includegraphics[height=35mm]{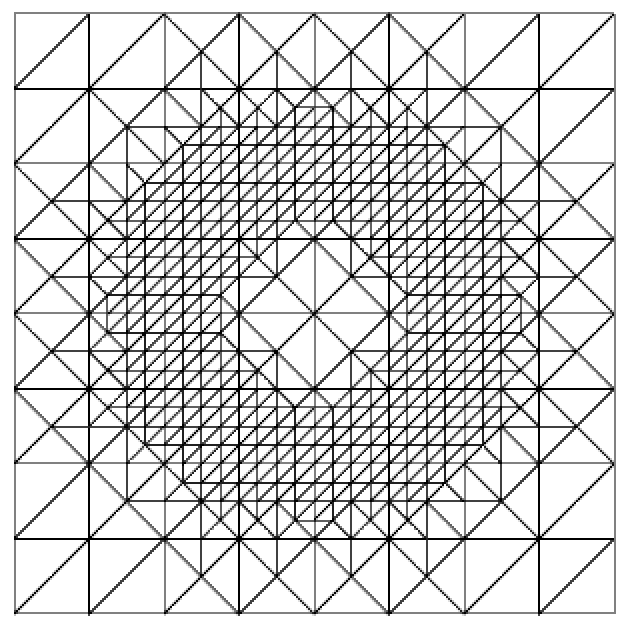} &
\includegraphics[height=35mm]{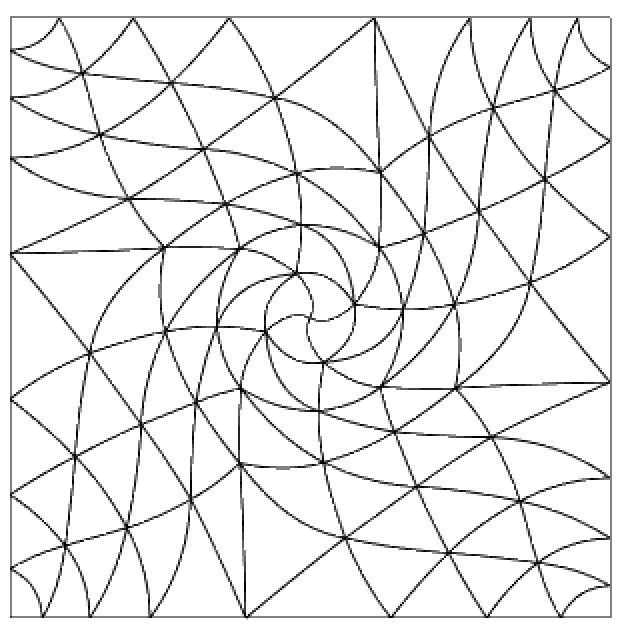} &
\includegraphics[height=35mm]{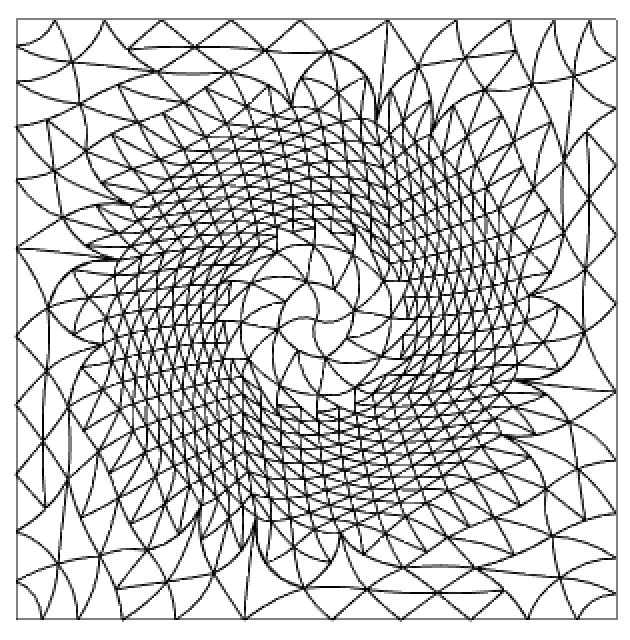}
\end{array}$
\end{center}
\vspace{-7mm}
\caption{Adaptivity on simplices: (left to right) (a) original mesh compared with the
optimized mesh using (b) \ha, (c) \ra, and (d) \hra.}
\label{fig_app_ex3}
\end{figure}

Figure \ref{fig_app_ex3} shows an example of a mesh ($N_E=128$) with simplices
optimized using \ha, \ra\ and \hra\ with $\mu_9^r$ and $\mu_{55}^h$ using the
targets in \eqref{eq_W_ex1}.  For this example, $F(\bx)$ is reduced 62.4\% by
\ha, 43.9\% by \ra, and 85.2\% by \hra, while the total element count is
increased to $N_E=928$ by \ha\ and to $N_E=1100$ by \hra.  Note that while the
final element count is significantly higher for the examples in Fig.
\ref{fig_app_ex1}-\ref{fig_app_ex3} for \hra\ in comparison to \ra, that is
expected due to the target size and aspect-ratio.  Additionally, as we have
demonstrated in Fig. \ref{fig_pois_error}, the total number of degrees of
freedom required for a desired accuracy in solution is typically lower for
meshes obtained with \hra\ as compared to \ra\ and \ha.

%==========================

\subsection{Mesh adaptivity with derefinement} \label{subsec_deref}
The examples in the previous sections demonstrate that the ability of
\emph{adding elements} with \hra\ can be critical for satisfying size and
aspect-ratio adaptivity targets, depending on the topology, size, and shape of
the elements in the original mesh.  For time-dependent problems however, where
the need for spatial resolution changes with time, the ability to \emph{remove
elements} (i.e, derefinement) is also critically important for maximizing the
computational efficiency of a calculation.

\begin{figure}[b!]
\begin{center}
$\begin{array}{ccc}
\includegraphics[height=35mm]{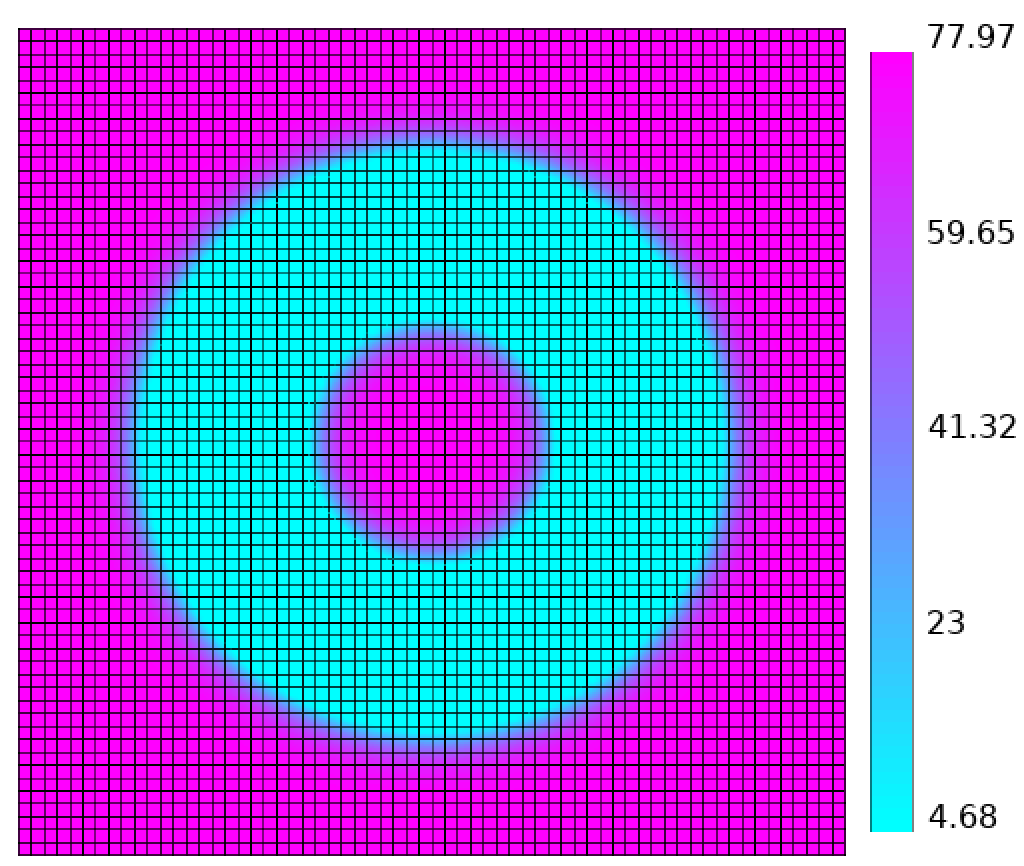} &
\includegraphics[height=35mm]{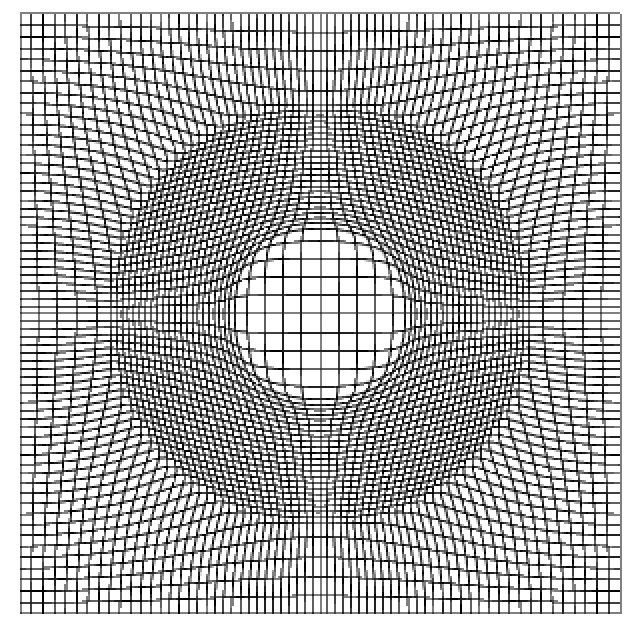} &
\includegraphics[height=35mm]{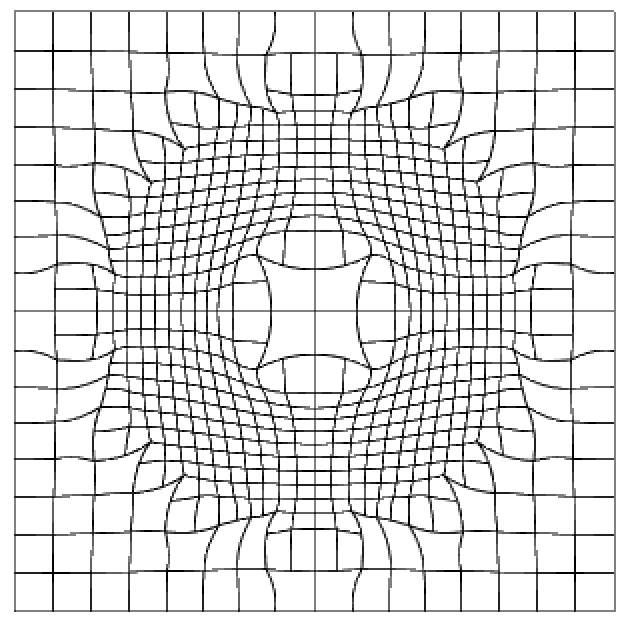}  \\
\end{array}$
\end{center}
\vspace{-7mm}
\caption{Adaptivity with derefinement: (left) original mesh compared with the
optimized mesh using (center) \ra and (right) \hra.}
\label{fig_app_ex4}
\end{figure}

Here, we consider the example from Section \ref{subsec_analytic} where we
optimize the mesh using $\mu_9^r$ and $\mu^h_{55}$ with the targets
\eqref{eq_W_ex1}.  We start with a $4 \times 4$ mesh and isotropically refine
each element four times, to obtain a mesh with $N_E=4096$ elements. As a
result, the final mesh has many more elements than required.  In contrast to
Fig. \ref{fig_app_ex1}, where the mesh was much coarser and the resulting
metric values were highest for the elements in the annulus region with elements
larger than specified by $W$, the metric values in Fig. \ref{fig_app_ex4} (left
panel) are highest for the elements outside the annulus region, as those
elements are much smaller than required.  With \ra, the element shape and size
are optimized to reduce the TMOP objective function by 55.4\% but the domain
still has more resolution than desired. In contrast, \hra\ removes the unneeded
resolution by reducing the total element count to $N_E=664$, which reduces
$F(\bx)$ by 98.6\%.

%==========================

\subsection{2D and 3D \hra\ in ALE hydrodynamics}

In this section, we consider a 2D and a 3D test case for \hra\ in the context
of Arbitrary Lagrangian-Eulerian (ALE) hydrodynamics.  The 2D case represents a
high velocity impact of gasses that was originally proposed in \cite{Barlow14}.
It involves three materials that represent an {\em impactor}, a {\em wall}, and
the {\em background}.  This problem is used to demonstrate the method's
behavior in an impact simulation that cannot be executed in Lagrangian frame to
final time as it produces large mesh deformations.  The complete thermodynamic
setup of this problem and additional details about our multi-material finite
element discretization and overall ALE method can be found in \cite{Dobrev2016,
Dobrev2018}. This method is available in a multi-material ALE code where we
have integrated the TMOP-based $r-$refinement framework in previous work
\cite{TMOP2020}.

To demonstrate the effectiveness of \hra\ in capturing the interface of
different material regions, we take snapshots of the discrete ALE solution from
the gas-impactor test and use these solutions to adapt the corresponding
meshes.  Following the approach described in Section \ref{subsec_pois}, we use
the gradient of the discrete density solution to set the size and aspect-ratio
targets and then optimize the mesh using $\mu^r_9$ and $\mu_1^h$ \shsz\
metrics.

\begin{figure}[b!]
\centerline{
  \includegraphics[height=0.3\textwidth]{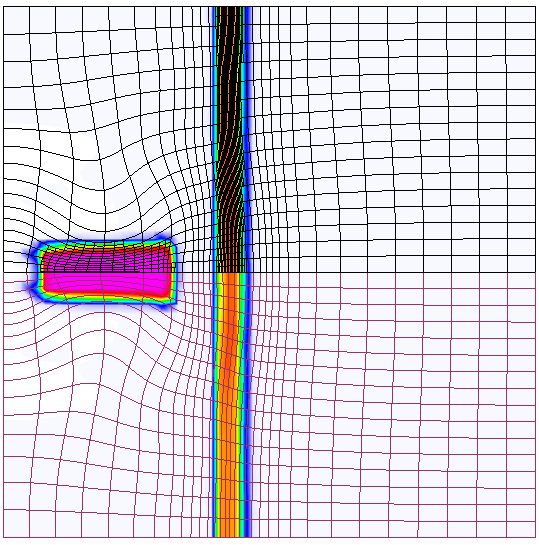} \hspace{0.3mm}
  \includegraphics[height=0.3\textwidth]{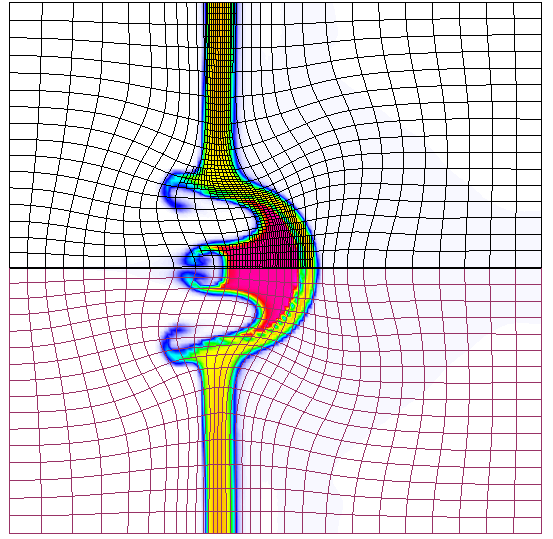} \hspace{0.3mm}
  \includegraphics[height=0.3\textwidth]{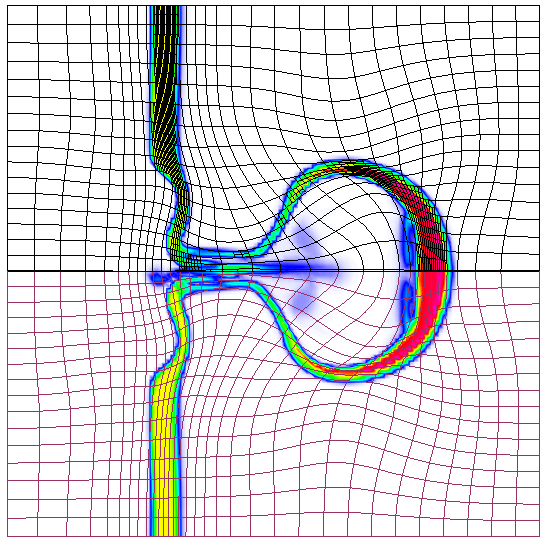}}
\caption{Mesh adaptivity in 2D ALE hydrodynamics: time evolution of density and mesh position at times
  0.5 (left), 4 (center) and 9 (right) in the 2D gas impact test case. The
  meshes obtained with \hra\ and \ra\ are shown in the upper (black) and lower (red) half of the
  snapshot, respectively. Notice the increased resolution in the {\em wall} and the {\em impactor}
  with \hra.}
\label{fig_gas_impact}
\end{figure}

Figure \ref{fig_gas_impact} shows a comparison of the optimized meshes obtained
using $r-$ and \hra\ along with the density function, before, during, and after
the initial impact. These results clearly show that we can achieve a much finer
resolution in the regions of interest where different materials are interacting
with \hra\ as compared to \ra.  This example also demonstrates the importance
of derefinement, as we can see that the regions that require isotropic and
anisotropic refinement are very different at different times.

\begin{figure}[b!]
\begin{center}
$\begin{array}{cc}
\includegraphics[height=0.3\textwidth]{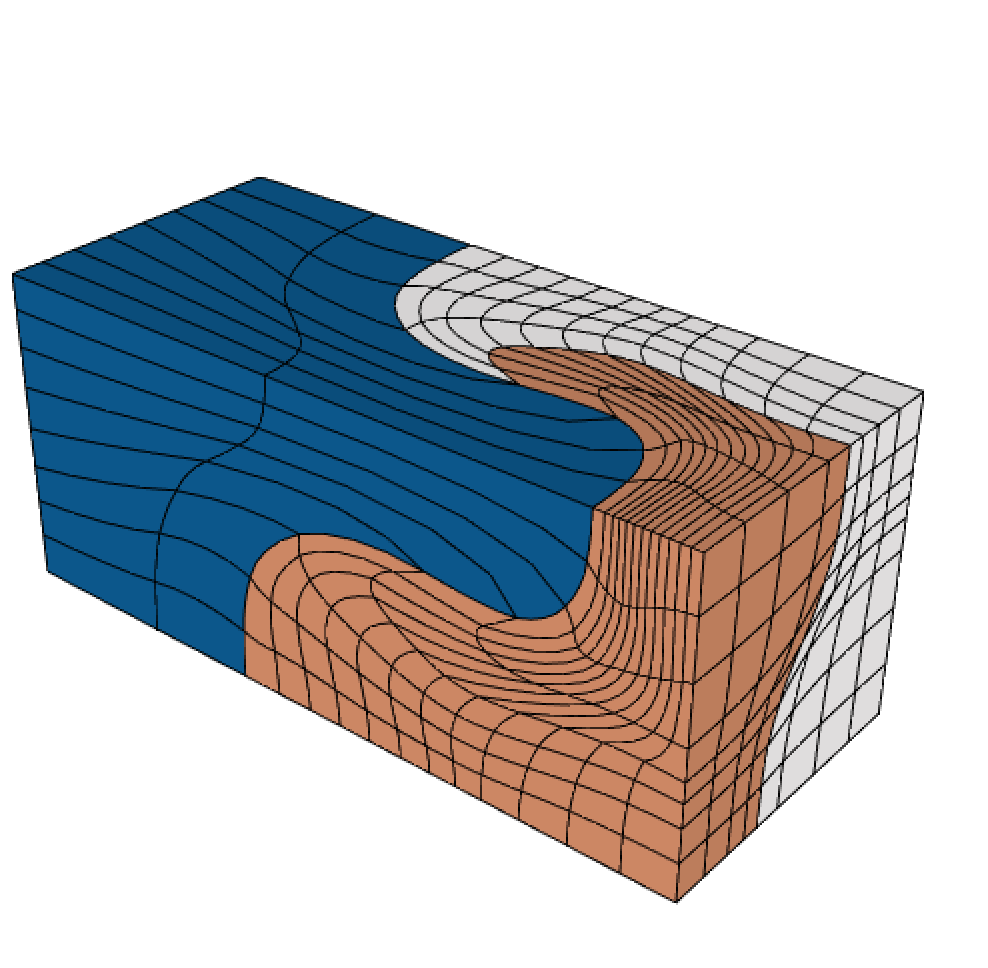} &
\includegraphics[height=0.3\textwidth]{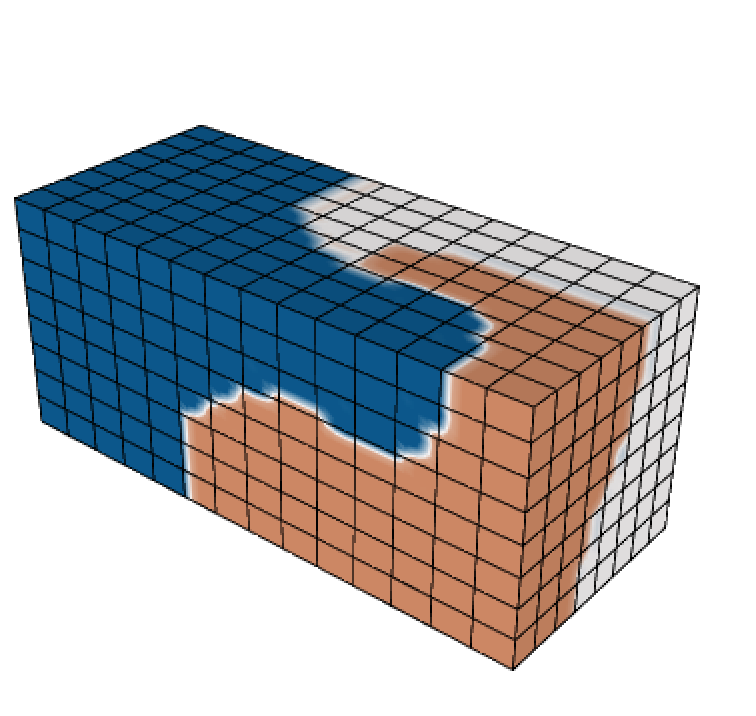} \\
\textrm{(a) Material indicator from Laghos.} & \textrm{(b) Uniform hexahedron mesh.} \\
\includegraphics[height=0.3\textwidth]{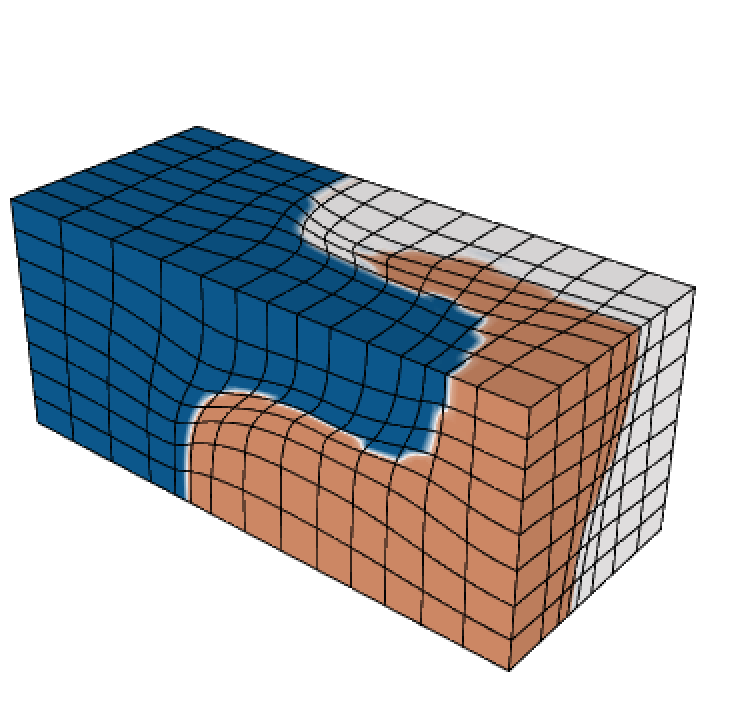} &
\includegraphics[height=0.3\textwidth]{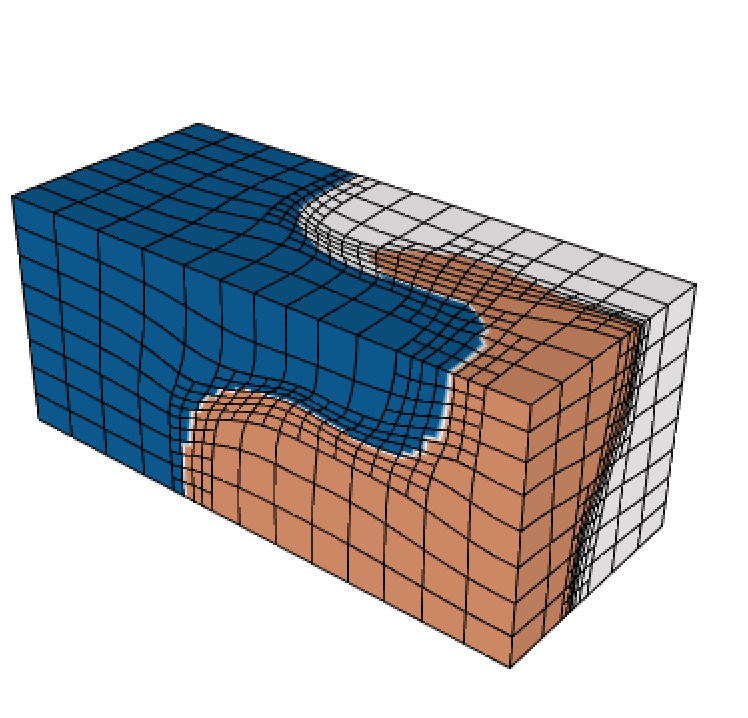} \\
\textrm{(c) Mesh optimized using \ra.} & \textrm{(d) Mesh optimized using \hra.} %% \\
%% \multicolumn{2}{c}{\includegraphics[height=0.075\textwidth]{figures/laghoslegend}}
\end{array}$
\end{center}
\vspace{-7mm}
\caption{Mesh adaptivity in 3D ALE hydrodynamics: (a) material indicator
function obtained from Laghos is mapped to a (b) uniform hexahedron mesh that
is used for mesh optimization.  A comparison of the meshes obtained using (c)
\ra\ and (d) \hra\ is shown below. Different colors indicate different
materials.}
\label{fig_triple_point}
\end{figure}

For the 3D test case, we consider a generalization of the triple point problem
in the Lagrangian High-Order Solver (Laghos) miniapp \cite{laghos}, where the
domain ($\Omega = [0, 7] \times [0, 3] \times [0, 3]$) is split into 5
different regions with 3 different materials. A complete thermodynamical setup
of this problem is discussed in Section 7.4 of \cite{zeng2016variational}, and
a snapshot of the material regions at simulation time $t=5$ is shown in Fig.
\ref{fig_triple_point}(a).  Since Laghos solves the time-dependent Euler
equations in a moving Lagrangian frame of reference, the initial mesh produces
no material diffusion.  However, in practical simulations, large mesh
deformations require ALE methods and mesh optimization. Thus, it is important
to be able to adapt the mesh in order to capture curved material interfaces.

For the purposes of this test, we start by transferring the material indicator
function defined on the Laghos mesh to a uniform hexahedron mesh (Fig.
\ref{fig_triple_point}(b)), using a high-order interpolation library
\cite{mittal2019nonconforming}.  The target shape is set to be the same as that
of an ideal element and the target size is computed using the magnitude of the
gradient of the material indicator. The \ra\ component is based on a \shsz\
metric ($\mu^r_{321}$) and the $h-$refinement decisions are made using a \sz\
metric ($\mu^h_{315}$).  Figure \ref{fig_triple_point} shows a comparison of
the meshes obtained from \ra\ and \hra. While \ra\ improves the original mesh
to align with the original material indicator, its effectiveness in reducing
the element size is limited due to the topology of the mesh. In contrast, \hra\
allows us to satisfy the size targets by refining the elements that are at the
interface of different materials. The uniform hexahedron mesh used here has
$N_E=896$ elements and \hra\ increases the total element count to 3724. To
quantify the comparison between these meshes, we look at how well the original
material indicator can be represented on each of the meshes in comparison to
the mesh from Laghos. Using high-order interpolation, we compute a volumetric
error: \begin{eqnarray} e = \int_{\Omega} \bigg(\eta_L(\bx_L) -
\eta(\bx_s)\bigg)^2, \end{eqnarray} that depends on the material indicator
$\eta_L$ defined on the mesh obtained from Laghos ($\bx_L$) and the same
material indicator function mapped on to the optimized mesh ($\bx_s$). The
volumetric error is 1.71 for the uniform hexahedral mesh, 0.62 for the \ra\
mesh, and 0.078 for the \hra\ mesh. Thus, we see that the volumetric error is
almost an order of magnitude lower for the \hra\ mesh in comparison to \ra\ and
about $20 \times$ lower in comparison to the uniform hexahedron mesh.

The examples presented in this section show the importance of \hra\ in the
context of complex Lagrangian hydrodynamics applications. In future work, we
plan to integrate the \hra\ framework in high-order applications with moving
meshes to demonstrate the runtime benefits in computational performance.

\section{Conclusion} \label{sec_conc}
In this paper we presented a new approach for \hra\ of high-order meshes. The
\ra\ component of our method is based on the Target Matrix Optimization
Paradigm, where the mesh is optimized to minimize a quality metric that depends
on the transformation from the current-to-target coordinates. For \ha, we
introduced a TMOP-based quality estimator that determines the refinement type
(isotropic or anisotropic) on an element-by-element basis to reduce the TMOP
objective function.  The \hra\ procedure uses an iterative loop of
$r-$refinements followed by $h-$refinements to adapt the mesh while minimizing
the computational cost of the calculation due to increase in the element count.
The new \hra\ algorithms are freely available in the MFEM finite element
library \cite{mfem-web}.  Numerical experiments show that \hra\ helps satisfy
the \emph{size} and \emph{aspect-ratio} targets where $h-$ and \ra\ cannot, due
to the topology of the input mesh. Using a 2D problem with a known exact
solution, we demonstrated that \hra\ requires many fewer degrees of freedom for
the desired accuracy, in comparison to $h-$ and \ra. We also used several
simple 2D problems to demonstrate \hra\ for analytic function-based adaptivity
in meshes with quadrilaterals and simplices, and the ability to derefine, which
is crucial for time-dependent problems where spatial resolution requirements
can change with time.  Finally, a 2D and a 3D Lagrangian hydrodynamics-based
test problems were used to compare the effectiveness of $r-$ with
$hr-$refinement in capturing complex curvilinear multi-material interfaces. In
the future, we plan to investigate methods for combined $r\!p-$ and
$hr\!p-$adaptivity.  We will also evaluate the current \hra\ framework in
applications with high-order moving meshes.

\noindent \section*{Disclaimer}
This document was prepared as an account of work sponsored by an agency of the
United States government. Neither the United States government nor Lawrence
Livermore National Security, LLC, nor any of their employees makes any
warranty, expressed or implied, or assumes any legal liability or
responsibility for the accuracy, completeness, or usefulness of any
information, apparatus, product, or process disclosed, or represents that its
use would not infringe privately owned rights. Reference herein to any specific
commercial product, process, or service by trade name, trademark, manufacturer,
or otherwise does not necessarily constitute or imply its endorsement,
recommendation, or favoring by the United States government or Lawrence
Livermore National Security, LLC. The views and opinions of authors expressed
herein do not necessarily state or reflect those of the United States
government or Lawrence Livermore National Security, LLC, and shall not be used
for advertising or product endorsement purposes.

\bibliographystyle{elsarticle-num}
\bibliography{hr}

\newcommand{\noopsort}[1]{} \newcommand{\printfirst}[2]{#1}
  \newcommand{\singleletter}[1]{#1} \newcommand{\switchargs}[2]{#2#1}
\begin{thebibliography}{10}
\expandafter\ifx\csname url\endcsname\relax
  \def\url#1{\texttt{#1}}\fi
\expandafter\ifx\csname urlprefix\endcsname\relax\def\urlprefix{URL }\fi
\expandafter\ifx\csname href\endcsname\relax
  \def\href#1#2{#2} \def\path#1{#1}\fi

\bibitem{TMOP2020}
V.~A. Dobrev, P.~Knupp, T.~V. Kolev, K.~Mittal, R.~N. Rieben, V.~Z. Tomov,
  Simulation-driven optimization of high-order meshes in {ALE} hydrodynamics,
  Comput. Fluids (2020).

\bibitem{vollmer1999improved}
J.~Vollmer, R.~Mencl, H.~Mueller, Improved laplacian smoothing of noisy surface
  meshes, in: Computer graphics forum, Vol.~18, Wiley Online Library, 1999, pp.
  131--138.

\bibitem{field1988laplacian}
D.~A. Field, Laplacian smoothing and delaunay triangulations, Communications in
  applied numerical methods 4~(6) (1988) 709--712.

\bibitem{taubin2001linear}
G.~Taubin, et~al., Linear anisotropic mesh filtering, Res. Rep. RC2213 IBM
  1~(4) (2001).

\bibitem{Knupp2012}
P.~Knupp, Introducing the target-matrix paradigm for mesh optimization by node
  movement, Engr. with Comptr. 28~(4) (2012) 419--429.

\bibitem{gargallo2015optimization}
A.~Gargallo-Peir{\'o}, X.~Roca, J.~Peraire, J.~Sarrate, Optimization of a
  regularized distortion measure to generate curved high-order unstructured
  tetrahedral meshes, International Journal for Numerical Methods in
  Engineering 103~(5) (2015) 342--363.

\bibitem{mittal2019mesh}
K.~Mittal, P.~Fischer, Mesh smoothing for the spectral element method, Journal
  of Scientific Computing 78~(2) (2019) 1152--1173.

\bibitem{dobrev2019target}
V.~Dobrev, P.~Knupp, T.~Kolev, K.~Mittal, V.~Tomov, The target-matrix
  optimization paradigm for high-order meshes, SIAM Journal on Scientific
  Computing 41~(1) (2019) B50--B68.

\bibitem{Greene2017}
P.~T. Greene, S.~P. Schofield, R.~Nourgaliev, Dynamic mesh adaptation for front
  evolution using discontinuous {G}alerkin based weighted condition number
  relaxation, Journal of Computational Physics 335 (2017) 664--687.

\bibitem{Peiro2018}
M.~Turner, J.~Peir\'{o}, D.~Moxey, Curvilinear mesh generation using a
  variational framework, Computer-Aided Design 103 (2018) 73--91.

\bibitem{Huang1994}
W.~Huang, Y.~Ren, R.~D. Russell, Moving mesh partial differential equations
  ({MMPDES}) based on the equidistribution principle, SIAM J. Numer. Anal.
  31~(3) (1994) 709--730.

\bibitem{Huang2010}
W.~Huang, R.~D. Russell, Adaptive moving mesh methods, Springer Science \&
  Business Media, 2010.

\bibitem{cerveny2019nonconforming}
J.~Cerveny, V.~Dobrev, T.~Kolev, Nonconforming mesh refinement for high-order
  finite elements, SIAM Journal on Scientific Computing 41~(4) (2019)
  C367--C392.

\bibitem{barros2004error}
F.~Barros, S.~Proen{\c{c}}a, C.~de~Barcellos, On error estimator and
  p-adaptivity in the generalized finite element method, International Journal
  for Numerical Methods in Engineering 60~(14) (2004) 2373--2398.

\bibitem{mackenzie2020hr}
J.~Mackenzie, W.~Mekwi, An hr-adaptive method for the cubic nonlinear
  schr{\"o}dinger equation, Journal of Computational and Applied Mathematics
  364 (2020) 112320.

\bibitem{piggott2005h}
M.~Piggott, C.~Pain, G.~Gorman, P.~Power, A.~Goddard, h, r, and hr adaptivity
  with applications in numerical ocean modelling, Ocean modelling 10~(1-2)
  (2005) 95--113.

\bibitem{jahandari2020finite}
H.~Jahandari, S.~MacLachlan, R.~D. Haynes, N.~Madden, Finite element modelling
  of geophysical electromagnetic data with goal-oriented hr-adaptivity,
  Computational Geosciences 24 (2020) 1257--1283.

\bibitem{piggott2009anisotropic}
M.~Piggott, P.~Farrell, C.~Wilson, G.~Gorman, C.~Pain, Anisotropic mesh
  adaptivity for multi-scale ocean modelling, Philosophical Transactions of the
  Royal Society A: Mathematical, Physical and Engineering Sciences 367~(1907)
  (2009) 4591--4611.

\bibitem{ong2013hr}
B.~Ong, R.~Russell, S.~Ruuth, An hr moving mesh method for one-dimensional
  time-dependent pdes, in: Proceedings of the 21st International Meshing
  Roundtable, Springer, 2013, pp. 39--54.

\bibitem{mostaghimi2015anisotropic}
P.~Mostaghimi, J.~R. Percival, D.~Pavlidis, R.~J. Ferrier, J.~L. Gomes, G.~J.
  Gorman, M.~D. Jackson, S.~J. Neethling, C.~C. Pain, Anisotropic mesh
  adaptivity and control volume finite element methods for numerical simulation
  of multiphase flow in porous media, Mathematical Geosciences 47~(4) (2015)
  417--440.

\bibitem{piggott2008new}
M.~Piggott, G.~Gorman, C.~Pain, P.~Allison, A.~Candy, B.~Martin, M.~Wells, A
  new computational framework for multi-scale ocean modelling based on adapting
  unstructured meshes, International Journal for Numerical Methods in Fluids
  56~(8) (2008) 1003--1015.

\bibitem{walkley2002anisotropic}
M.~Walkley, P.~K. Jimack, M.~Berzins, Anisotropic adaptivity for the finite
  element solutions of three-dimensional convection-dominated problems,
  International journal for numerical methods in fluids 40~(3-4) (2002)
  551--559.

\bibitem{antonietti2020hr}
P.~F. Antonietti, P.~Houston, An hr-adaptive discontinuous galerkin method for
  advection-diffusion problems, in: Communications to SIMAI congress, Vol.~3,
  2020.

\bibitem{edwards1993h}
M.~G. Edwards, J.~T. Oden, L.~Demkowicz, An h-r-adaptive approximate riemann
  solver for the euler equations in two dimensions, SIAM Journal on scientific
  computing 14~(1) (1993) 185--217.

\bibitem{knupp2019target}
P.~Knupp, Target formulation and construction in mesh quality improvement,
  Tech. Rep. LLNL-TR-795097, Lawrence Livermore National Lab.(LLNL), Livermore,
  CA (United States) (2019).

\bibitem{IMR2018}
V.~A. Dobrev, P.~Knupp, T.~V. Kolev, V.~Z. Tomov, Towards Simulation-Driven
  Optimization of High-Order Meshes by the {T}arget-{M}atrix {O}ptimization
  {P}aradigm, Springer International Publishing, 2019, pp. 285--302.

\bibitem{babuska1979reliable}
I.~Babuska, W.~C. Rheinboldt, Reliable error estimation and mesh adaptation for
  the finite element method., Tech. rep., Maryland Univ College Park Inst For
  Physical Science And Technology (1979).

\bibitem{ainsworth2011posteriori}
M.~Ainsworth, J.~T. Oden, A posteriori error estimation in finite element
  analysis, Vol.~37, John Wiley \& Sons, 2011.

\bibitem{mittal2019nonconforming}
K.~Mittal, S.~Dutta, P.~Fischer, Nonconforming {Schwarz}-spectral element
  methods for incompressible flow, Computers \& Fluids 191 (2019) 104237.

\bibitem{mfem}
R.~Anderson, J.~Andrej, A.~Barker, J.~Bramwell, J.-S. Camier, J.~C.~V. Dobrev,
  Y.~Dudouit, A.~Fisher, T.~Kolev, W.~Pazner, M.~Stowell, V.~Tomov,
  I.~Akkerman, J.~Dahm, D.~Medina, S.~Zampini, {MFEM}: A modular finite element
  library, Computers \& Mathematics with Applications (2020).
\newblock \href {https://doi.org/10.1016/j.camwa.2020.06.009}
  {\path{doi:10.1016/j.camwa.2020.06.009}}.

\bibitem{mfem-web}
{MFEM}: Modular finite element methods {[Software]}, \url{https://mfem.org}.
\newblock \href {https://doi.org/10.11578/dc.20171025.1248}
  {\path{doi:10.11578/dc.20171025.1248}}.

\bibitem{Barlow14}
A.~Barlow, R.~Hill, M.~J. Shashkov, Constrained optimization framework for
  interface-aware sub-scale dynamics closure model for multimaterial cells in
  {L}agrangian and arbitrary {L}agrangian-{E}ulerian hydrodynamics, J. Comput.
  Phys. 276~(0) (2014) 92--135.

\bibitem{Dobrev2016}
V.~A. Dobrev, T.~V. Kolev, R.~N. Rieben, V.~Z. Tomov, Multi-material closure
  model for high-order finite element {L}agrangian hydrodynamics, Int. J.
  Numer. Meth. Fluids 82~(10) (2016) 689--706.

\bibitem{Dobrev2018}
R.~W. Anderson, V.~A. Dobrev, T.~V. Kolev, R.~N. Rieben, V.~Z. Tomov,
  High-order multi-material {ALE} hydrodynamics, SIAM J. Sci. Comp. 40~(1)
  (2018) B32--B58.

\bibitem{laghos}
{Laghos}: High-order {L}agrangian hydrodynamics miniapp {[Software]},
  \url{https://github.com/ceed/Laghos} (2020).

\bibitem{zeng2016variational}
X.~Zeng, G.~Scovazzi, A variational multiscale finite element method for
  monolithic ale computations of shock hydrodynamics using nodal elements,
  Journal of Computational Physics 315 (2016) 577--608.

\end{thebibliography}

\end{document}